\documentclass[a4paper,10pt]{amsart}
\usepackage[latin1]{inputenc}
\usepackage{amsmath}
\usepackage{amsfonts}
\usepackage{amssymb}
\usepackage[T1]{fontenc}
\usepackage{dsfont}
\usepackage{float}
\usepackage{lscape}
\usepackage[all]{xy}
\usepackage{stmaryrd}

\newtheorem{theorem}{Theorem}[section]
\newtheorem{prop}[theorem]{Proposition}
\newtheorem{lem}[theorem]{Lemma}
\newtheorem{corol}[theorem]{Corollary}

\theoremstyle{definition}\newtheorem{defi}[theorem]{Definition}

\newtheorem{exmp}[theorem]{Example}
\newtheorem{rmq}[theorem]{Remark}

\def\N{{\mathbb{N}}}
\def\Z{{\mathbb{Z}}}
\def\C{{\mathbb{C}}}

\def\Q{{\mathbb{Q}}}

\def\A{{\mathbb{A}}}
\def\B{{\mathbb{B}}}
\def\D{{\mathbb{D}}}
\def\E{{\mathbb{E}}}

\def\si{{\rm{\, if\, }}}
\def\ou{{\rm{\, or\,  }}}
\def\et{{\rm{\, and\, }}}
\def\sinon{{\rm{\, otherwise\, }}}

\def\ssi{ if and only if }

\def\rep{{\rm{rep}}}
\def\CC{{\mathcal{C}}}

\def\ddim{{\textbf{dim}}\,}

\def\ind{{\textrm{ind}}}
\def\modd{{\textrm{-mod}}}
\def\Ob{{\rm{Ob}}}

\def\Hom{{\rm{Hom}}}
\def\Ext{{\rm{Ext}}}

\def\stab{{\rm{stab}}}

\def\nsl{non-simply-laced }
\def\sl{simply-laced }

\def\Cl{{\rm{Cl}}}

\def\Mut{{\rm{Mut}}}

\def\barQ{{\overline Q}}

\def\diag{{\rm{diag}}}
\def\ens#1{\left\{ #1 \right\}}
\def\fl{{\longrightarrow}\,}

\def\<{\left<}
\def\>{\right>}

\def\Gr{{\rm{Gr}}}

\def\Daffine{\tilde{\mathbb D}}

\title{An approach to non-simply-laced cluster algebras}
\author{\textsc{G. Dupont}}
 \address{Universit\'e de Lyon \\
Universit\'e Lyon 1 \\
Institut Camille Jordan CNRS UMR 5208 \\
43, boulevard du 11 novembre 1918\\
F-69622 Villeurbanne Cedex.}
 \email{dupont@math.univ-lyon1.fr}

\begin{document}
\maketitle

\begin{abstract}
	Let $\Delta$ be an oriented valued graph equipped with a group of admissible automorphisms satisfying a certain stability condition. We prove that the (coefficient-free) cluster algebra $\mathcal A(\Delta/G)$ associated to the valued quotient graph $\Delta/G$ is a subalgebra of a quotient $\pi(\mathcal A(\Delta))$ of the cluster algebra associated to $\Delta$ by the action of $G$. When $\Delta$ is a Dynkin diagram, we prove that $\mathcal A(\Delta/G)$ and $\pi(\mathcal A(\Delta))$ coincide. As an example of application, we prove that affine valued graphs are mutation-finite, giving an alternative proof to a result of Seven.
\end{abstract}

\setcounter{tocdepth}{2}
\tableofcontents

\section*{Introduction}
	Cluster algebras were introduced in 2001 by Fomin and Zelevinsky in order to define a combinatorial framework for the study of total positivity in algebraic groups and canonical bases in quantum groups \cite{cluster1}. Since they have found applications in various areas of mathematics like combinatorics, Teichm\"uller theory, Lie theory, Poisson geometry or representation theory. 
	
	Most of the advanced results in the theory of cluster algebras come from categorifications using cluster categories of acyclic quivers or preprojective algebras \cite{BMRRT, CC, CK2, GLS:rigid2}. Nevertheless, these categorifications are both mainly restricted to the case where the cluster algebra is simply-laced.

	If one wants to find similar results on non-simply-laced cluster algebras, a first approach would consist in using representation theory of species (see e.g. \cite{Zhu:BGPcluster,Zhu:applications,zhu,hubery:cluster}). If this approach seems to be fully functional, it nevertheless asks to adapt most of the proofs and does not enable to use directly the numerous existing results about simply-laced cluster algebras.
	
	In Kac-Moody-Lie  theory, algebras associated to non-simply-laced diagrams are often studied as quotients of algebras equipped with a group of automorphisms (see \cite{Humphreys:Lie} for example). In general, the study of non-simply-laced valued datas using simply-laced valued datas with automorphisms was developed bu Lusztig and Kac \cite{Lusztig:quantumbook,kac}. This construction is also known to be fruitful at the level of representations of quivers \cite{hubery}. One of the main interests of this method is that it allows to transport at the level of \nsl datas results concerning \sl datas.
	
	In this paper, we use this construction of \nsl diagrams in order to obtain results at the level of \nsl cluster algebras using those concerning \sl cluster algebras. Thus, if $\Delta$ is a valued graph equipped with a group of admissible automorphisms $G$ satisfying a certain stability condition, we prove that the cluster algebra associated to the quotient graph $\Delta/G$ is a subalgebra of a certain quotient of the cluster algebra $\mathcal A(\Delta)$ with respect to an action of $G$ on the ring of Laurent polynomials containing the cluster algebra. 
	
	Another approach, using Frobenius morphisms on species over finite fields was also developed by Yang \cite{yang}. It follows from results of Deng et Du that this method is actually related to the one we present here \cite{DD,DD2}. 
	
	Note also that a general categorification of non-simply-laced cluster algebras using quivers with automorphisms was recently developed by Demonet in the context of preprojective algebras \cite{Demonet, Demonet:PhD}.
	
	Figure \ref{figure:approchesnsl} sums up the situation
	\begin{figure}[H]
		\begin{picture}(260,120)(30,0)
			\put(0,80){$B \in M_q(\Z)$ skew-symmetric,}
			\put(0,70){$G \leq \mathfrak S_q$ group}
			\put(250,80){$\mathcal A(B)$}
			\put(250,70){simply-laced}
			
			\put(25,65){\vector(0,-1){35}}
			\put(30,50){\cite{kac,Lusztig:quantumbook}}
			
			\put(275,65){\vector(0,-1){35}}
			\put(280,50){\cite{Demonet,Demonet:PhD}}
			
			\put(160,80){\vector(1,0){80}}
			
			\put(160,90){\cite{BMRRT,CCS1,BMRT,CK1,CK2,GLS}$\ldots$}
			
			\put(160,20){\vector(1,0){80}}
			\put(160,5){\cite{Zhu:BGPcluster,Zhu:applications,zhu,hubery:cluster}$\ldots$}
			
			\put(0,20){$B/G \in M_s(\Z)$ skew-symmetrizable}
			\put(250,20){$\mathcal A(B/G)$}
			\put(250,10){non-simply-laced}
		\end{picture} 
		\caption{From \sl to \nsl cluster algebras}\label{figure:approchesnsl}
	\end{figure}
	
	The paper is organized as follows. Section \ref{section:preliminaries} recalls all the necessary definitions and background. In section \ref{section:mutationsorbitales}, we prove the main results of the paper (namely, Theorem \ref{theorem:mutationquotient} and Corollary \ref{corol:corolquotient}) with purely combinatorial methods. In section \ref{section:interpretationcategorique}, we give an interpretation in terms of invariant objects in cluster categories whenever it is possible. Finally, in section \ref{section:mutationfinitude}, we give an example of application to the problem of classification of valued graphs having a finite mutation class.
	
\begin{section}{Preliminaries}\label{section:preliminaries}
	\begin{subsection}{Cluster algebras}	
		A cluster algebra is defined from a skew-symmetrizable matrix, that is, a matrix $B \in M_q(\Z)$ such that there exists a diagonal matrix $D \in M_q(\Z)$ with non-negative entries such that $DB$ is skew-symmetric. A \emph{seed} is a pair $(B,\textbf u)$ where $B=(b_{ij}) \in M_q(\Z)$ is a skew-symmetrizable matrix and $\textbf u=\ens{u_1, \ldots, u_q}$ is a set of indeterminates over $\Q$. The matrix $B$ is called the \emph{exchange matrix} of the seed $(B,\textbf u)$ and $\textbf u$ is called the \emph{cluster} of the seed $(B,\textbf u)$.
		
		Given a seed $(B,\textbf u)$ and an index $k \in \ens{1, \ldots, q}$, the \emph{mutation of the seed in the direction $k$} is the new seed $\mu_k(B,\textbf u)=(B',\textbf u')$ where $B'=(b_{ij}') \in M_q(\Z)$ is given by
		$$b_{ij}'=\left\{\begin{array}{ll}
			-b_{ij} & \textrm{ if } i=k \textrm{ or } j=k,\\
			b_{ij}+\frac 12 \left(b_{ik}|b_{kj}|+|b_{ik}|b_{kj}\right) & \textrm{ otherwise.}
		\end{array}\right.$$
		and $\textbf u'=\textbf u \setminus \ens{u_k} \sqcup \ens{u_k'}$
		where $u_k'$ is related to $u_k$ by the following \emph{exchange relation}:
		$$u_ku_k'=\prod_{b_{ik}>0}u_i^{b_{ik}} + \prod_{b_{ik}<0}u_i^{-b_{ik}}.$$
		
		We denote by $\Mut(B,\textbf u)$ the set, called \emph{mutation class} of $(B,\textbf u)$, of all seeds that can be obtained from $(B,\textbf u)$ after a finite number of mutations. The seeds in $\Mut(B,\textbf u)$ are called \emph{mutation-equivalent} to $(B,\textbf u)$. We denote by $\Mut(B)$ the set, called \emph{mutation class} of $B$, of exchange matrices of seeds in $\Mut(B,\textbf u)$.
		
		The (coefficient-free) \emph{cluster algebra} with initial seed $(B,\textbf u)$ is the $\Z$-subalgebra $\mathcal A(B,\textbf u)$ of $\mathcal F=\Q(u_1, \ldots, u_q)$ generated by variables in clusters occurring in seeds mutation-equivalent to $(B,\textbf u)$. In other words,
		$$\mathcal A(B,\textbf u)=\Z[x \ : \ x \in \textbf c \textrm{ s.t. } \exists C \in M_q(\Z) \textrm{ s.t. } (C,\textbf c) \in \Mut(B,\textbf u)] \subset \mathcal F.$$
		
		The clusters occurring in seeds in $\Mut(B,\textbf u)$ are called the \emph{clusters} of the cluster algebra $\mathcal A(B,\textbf u)$ and $\textbf u$ is called the \emph{initial cluster}. The variables occurring in the clusters are called \emph{cluster variables}. We denote by $\Cl(B,\textbf u)$ the set of cluster variables in $\mathcal A(B,\textbf u)$.
		
		Fomin and Zelevinsky proved that $\mathcal A(B,\textbf u)$ is a subring of the ring of Laurent polynomials in $\textbf u$. This is referred to as the \emph{Laurent phenomenon}. Thus, given an element $x \in \mathcal A(B,\textbf u)$, we can write it
		$$x=\frac{P(u_1, \ldots, u_q)}{u_1^{n_1} \ldots u_q^{n_q}}$$
		such that $P(u_1, \ldots, u_q)$ is a polynomial not divisible by any of the $u_i$ and $n_i \in \Z$ for every $i \in \ens{1, \ldots, q}$. The sequence $\delta(x)=(n_1, \ldots, n_q)$ is called the \emph{denominator vector of $x$}.
		
		When there is no confusion, we will omit the reference to the initial cluster in the notations. 
		
		If a matrix $B \in M_q(\Z)$ is skew-symmetric, then all the matrices in $\Mut(B)$ are also skew-symmetric. In this case, the cluster algebra $\mathcal A(B)$ is called \emph{simply-laced}. If $B$ is skew-symmetrizable but not skew-symmetric, $\mathcal A(B)$ is called a \emph{non-simply-laced}.
		
		Given a skew-symmetrizable matrix $B$, we can define the \emph{Cartan counterpart $C(B)=(c_{ij})$ of $B$} by setting $c_{ii}=2$ and $c_{ij}=-|b_{ij}|$. This defines a generalized Cartan matrix to which we can associate a valued diagram $\Gamma(C(B))$ (see \cite{kac} for example). The valued graph $Q_B$ associated to $B$ is thus the orientation of $\Gamma(C(B))$ given by $i \fl j$ if $b_{ij}>0$ for any $i,j \in \ens{1, \ldots, q}$. Thus, if $Q_B$ is the valued graph associated to $B$, we will sometimes write $\mathcal A(Q_B)$ instead of $\mathcal A(B)$.
	\end{subsection}

	\begin{subsection}{Cluster categories}
		If $B$ is skew-symmetric, then $Q=Q_B$ is a quiver with no loops and no 2-cycles. In this case, there is a fruitful framework for categorifying the cluster algebra $\mathcal A(Q)$ using the category of representations of $Q$. 
		
		More precisely, let $k \simeq \C$ be the field of complex numbers and $\rep(Q)$ be the category of $k$-representations of $Q$. We denote by $Q_0$ the set of vertices of $Q$ and $Q_1$ the set of arrows of $Q$. An element in $\rep(Q)$ is thus a pair $((V(i))_{i \in Q_0}, (V(\alpha))_{\alpha \in Q_1})$ such that each $V(i)$ is a finite-dimensional $k$-vector space and $V(\alpha):V(i) \fl V(j)$ is a $k$-linear map for every $\alpha: i \fl j$ in $Q_1$. As usual, we will identify $\rep(Q)$ with the category $kQ$-mod of finite-dimensional modules over the path algebra of $Q$. For every vertex $i \in Q_0$, we denote by $S_i$ the simple module associated to the vertex $i$, by $P_i$ its projective cover and by $I_i$ its injective hull.
		
		We denote by $D^b(kQ)$ the bounded derived category of $kQ$-mod with shift functor $[1]$ and Auslander-Reiten translation $\tau$. The functor $F=\tau^{-1}[1]$ is an auto-equivalence of $D^b(kQ)$. In \cite{BMRRT}, the authors defined the \emph{cluster category $\CC_Q$ of $Q$} as the orbit category of $F$ in $D^b(kQ)$. The objects in $\CC_Q$ are thus the objects in $D^b(kQ)$ and given two objects $M,N$ in $\CC_Q$, the morphisms from $M$ to $N$ in $\CC_Q$ are given by:
		$$\Hom_{\CC_Q}(M,N) = \bigoplus_{i \in \Z}\Hom_{D^b(kQ)}(M,F^iN).$$
		
		A \emph{module} in $\CC_Q$ is the image of a module under the composition of functors 
		$$kQ\modd \fl D^b(kQ) \fl \CC_Q$$
		where the first map sends a module $M$ to the corresponding complex concentrated in degree 0 and the second functors is the canonical map. We still denote by $M$ the image of a $kQ$-module $M$ in $\CC_Q$.
		
		The cluster category $\CC_Q$ is triangulated \cite{K}. Moreover, it is proved in \cite{BMRRT} that $\CC_Q$ is a Krull-Schmidt category whose isoclasses of indecomposable objects are given by
		$$\ind(\CC_Q)=\ind(kQ\modd) \sqcup \ens{P_i[1] \ : \ i \in Q_0}.$$
		Moreover, there is a duality
		$$\Ext^1_{\CC_Q}(M,N) \simeq D\Ext^1_{\CC_Q}(N,M)$$
		for any two objects $M,N$ in $\CC_Q$. This means that $\CC_Q$ is a 2-Calabi-Yau category.
		
		An object $M$ in $\CC_Q$ such that $\Ext^1_{\CC_Q}(M,M)=0$ is called \emph{rigid}. It is called a \emph{cluster-tilting object} if it is maximal rigid and if all its distinct indecomposable direct summands are pairwise non-isomorphic. 
	\end{subsection}

	\begin{subsection}{The Caldero-Chapoton map}
		The Caldero-Chapoton map is a map from the set of objects in the cluster category $\CC_Q$ taking its value in the ring of Laurent polynomials in $\textbf u$ containing the cluster algebra $\mathcal A(Q,\textbf u)$. Among many interests, this map allows to realize explicitly the categorification of $\mathcal A(Q)$ using $\CC_Q$. Before introducing it precisely, we will need some more notations.
		
	 	We denote by $K_0(kQ)$ the Grothendieck group of $kQ$-mod, that is, the free abelian group over the isoclasses of $kQ$-modules with relations $M+N=X$ for any short exact sequence $0 \fl M \fl X \fl N \fl 0$. The \emph{dimension vector} of a representation $M$ of $Q$ is the vector
		$$\ddim M=\left( \dim M(i)\right)_{i \in Q_0} \in \N^{Q_0}.$$
		We will denote by $\alpha_i$ the $i$-th vector of the canonical basis of $\Z^{Q_0}$. Thus, for every $i \in Q_0$, $\ddim S_i=\alpha_i$ and the map $\ddim$ induces an isomorphism of abelian groups 
		$$\ddim : K_0(kQ) \xrightarrow{\sim} \Z^{Q_0}.$$
		
		As $kQ$ is hereditary, the Euler form on $kQ$-mod is given by
		$$\<M,N\>=\dim \Hom_{kQ}(M,N)-\dim \Ext^1_{kQ}(M,N)$$
		for any $kQ$-modules $M,N$. It is well defined on the Grothendieck group.
		
		Given a $kQ$-module $M$ and a dimension vector $\textbf e$, we denote by
		$$\Gr_{\textbf e}(M)=\ens{N \textrm{ is a submodule of }M \ : \ \ddim N = \textbf e}$$ 
		the grassmannian of submodules of $M$ of dimension $\textbf e$. This is a closed subset of the usual grassmannian of $k$-vector-spaces. It it thus a projective variety and we can consider its Euler-Poincar\'e characteristic $\chi(\Gr_{\textbf e}(M))$.
		
		Roughly speaking, the Caldero-Chapoton map evaluated at a $kQ$-module $M$ is a generating series for Euler characteristics of grassmannians of submodules of $M$. More precisely, the definition is the following:
		\begin{defi}[\cite{CC}]
			The \emph{Caldero-Chapoton} map is the map $X_?:\Ob(\CC_Q) \fl \Z[\textbf u^{\pm 1}]$ given by:
			\begin{itemize}
			\item If $M,N$ are in $\Ob(\CC_Q)$, then $X_{M \oplus N}=X_MX_N;$
			\item If $M \simeq P_i[1]$, then $X_{P_i[1]}=u_i;$
			\item If $M$ is an indecomposable module, then
				\begin{equation}\label{eq:XM}
					X_M=\sum_{\textbf e \in K_0(kQ)} \chi(\Gr_{\textbf e}(M))\prod_{i \in Q_0} u_i^{-\<\textbf e, \ddim S_i\>-\<\ddim S_i, \ddim M - \textbf e\>}.
				\end{equation}
			\end{itemize}
		\end{defi}
		Note that equality (\ref{eq:XM}) holds for any $kQ$-module.
		
		The Caldero-Chapoton map realizes the categorification of $\mathcal A(Q)$ in the following sense:
		\begin{theorem}[\cite{CK2}]\label{theorem:correspondanceCK2}
			Let $Q$ be an acyclic quiver. Then $X_?$ induces a 1-1 correspondence
			$$\ens{\textrm{indecomposable rigid objects in }\CC_Q} \xrightarrow{\sim} \Cl(Q).$$
			Moreover, the map
			$$\left\{\begin{array}{rcl}
				\ens{\textrm{cluster-tilting objects in }\CC_Q} & \xrightarrow{\sim} & \ens{\textrm{clusters in }\mathcal A(Q)}\\
				T=\bigoplus_{i \in Q_0} T_i & \mapsto & \ens{X_{T_i} \ : \ i \in Q_0}
			\end{array}\right.$$
			is a 1-1 correspondence.
		\end{theorem}
	\end{subsection}
\end{section}

\begin{section}{Orbit mutations}\label{section:mutationsorbitales}
	In this section, we use Lusztig's approach to \nsl diagrams using \sl diagrams equipped with a group of automorphisms \cite{Lusztig:quantumbook}. We prove that under certain admissibility and stability assumptions, a \nsl cluster algebra can be realized as a subalgebra of a quotient of a \sl cluster algebra. Moreover, our methods allow to establish a link between \sl and \nsl cluster algebras of different types.
	
	\begin{subsection}{Automorphisms and quotient matrices}\label{section:automorphisms}
		Let $Q=(Q_0,Q_1)$ be a valued graph without $p$-cycles for $p \leq 2$. Let $B_{Q}=(b_{ij}) \in M_{Q_0}(\Z)$ be the corresponding skew-symmetrizable matrix. 
		\begin{defi}
			An elements $g \in \mathfrak S_{Q_0}$ is called an \emph{automorphism of $Q$} or \emph{automorphism of $B_{Q}$} if 
			$$b_{gi,gj}=b_{ij}$$
			for every $i,j \in Q_0$.
			
			A subgroup $G \leq \mathfrak S_{Q_0}$ is called \emph{an automorphism group of $B_Q$} if every $g \in G$ is an automorphism of $B_{Q}$.
			
			An automorphism group $G$ of $B_Q$ is called  \emph{admissible} if for every $i,j$ in the same $G$-orbit, there is no path from $i$ to $j$ of length $l \leq 2$. The pair $(B_Q, G)$ is then called an \emph{admissible pair}.
		\end{defi}
				
		Given a matrix $B \in M_{Q_0}(\Z)$ and an automorphism group $G$ of $B$, we denote by $\textbf i$ the $G$-orbit of an element $i \in Q_0$ and $\barQ_0$ the set
		$$\barQ_0=\ens{\textbf i \ : \ i \in Q_0}.$$
		of $G$-orbits in $Q_0$.
	
		\begin{exmp}\label{exmp:A3quotient}
			Let $Q$ be the Dynkin quiver of type $\A_3$ 
			$$\xymatrix{Q: 1 \ar@/_10pt/@{--}[rr]& \ar[l] 2 \ar[r] & 3}$$
			with associated matrix
			$$B=B_Q=\left[\begin{array}{rrr}
				0 & -1 & 0 \\
				1 & 0 & 1\\
				0 & -1 & 0 \\
			\end{array}\right]$$
			Let $G=\<(1,3)\>$ be the subgroup of $\mathfrak S_{\ens{1,2,3}}$ generated by the transposition $(1,3)$. Then $(B,G)$ is an admissible pair. In general, we will denote by dashed lines the action of the automorphism group on vertices of the quiver.
		\end{exmp}
		
		Given an admissible pair $(B,G)$, we define a matrix $B/G$, called \emph{quotient matrix} as follows:
		\begin{defi}
			Let $(A,G)$ be an admissible pair with $A=(a_{i,j}) \in M_{Q_0}(\Z)$, the \emph{quotient matrix} $A/G=(b_{\textbf i,\textbf j}) \in M_{\barQ_0}(\Z)$ is given by 
			$$b_{\textbf i, \textbf j}=\sum_{k \in \textbf i}a_{k,j}$$
			for every $\textbf i, \textbf j \in \barQ_0$.
		\end{defi}
		The matrix $A/G$ is well defined. Indeed, $b_{\textbf i,\textbf j}$ does not depend on the choice of $j \in \textbf j$ since $G$ is an automorphism group for $A$. An alternative way for describing the coefficients of the quotient matrix is the following:
		\begin{lem}
			Let $(A,G)$ be an admissible pair and $A/G=(b_{\textbf i, \textbf j})$. Then for every $\textbf i,\textbf j \in \barQ_0$, we have
			$$b_{\textbf i, \textbf j}=\frac{1}{|\stab_G(i)|}\sum_{g \in G} a_{gi,j}$$
			where $\stab_G(i)$ is the stabilizer of $i$ for the $G$-action on $Q_0$.
		\end{lem}
		\begin{proof}
			For every $g \in \stab(i)$, $gi=i$ and thus $a_{gi,j}=a_{i,j}$.
			$$\sum_{g \in G} a_{gi,j}=\sum_{g \in G/\stab(i)} |\stab(gi)| a_{gi,j}$$
			But $g \fl gi$ induces a 1-1 correspondence between $G/\stab(i)$ and $\textbf i$. Moreover, $\stab(i) \simeq \stab(gi)$ and thus
			$$\sum_{g \in G} a_{gi,j}= |\stab(i)| \sum_{k \in \textbf i} a_{k,j}$$
		\end{proof}
		
		We now prove that a quotient matrix is skew-symmetrizable:
		\begin{lem}\label{lem:symmetrizable}
			Let $(A,G)$ be an admissible pair. Then, $A/G$ is skew-symmetrizable.
		\end{lem}
		\begin{proof}
			Assume that $A \in M_{Q_0}(\Z)$ and set $B=A/G$. We shall find a diagonal matrix $\Delta \in M_{\barQ_0}(\Z)$ with non-negative coefficients such that $\Delta \overline B$ is skew-symmetric.
			
			As $A$ is skew-symmetrizable, there exists integers $(d_i)_{i \in Q_0}$ such that for every $i,j$, $d_ia_{ij}=-d_ja_{ji}$. Moreover, $G$ is an automorphism group of $A$, we can thus assume that $d_{gi}=d_i$ for every $i \in Q_0$ and we set $d_{\textbf i}=d_i$ if $\textbf i$ denotes the $G$-orbit of $i$.
			
			For every $\textbf i \in \barQ_0$, we set $\delta_{\textbf i}=d_{\textbf i}|\stab_G(i)|$ for an arbitrary vertex $i \in \textbf i$. Let $\Delta$ be the diagonal matrix $\diag(\delta_{\textbf i}, \textbf i \in \bar Q_0)$. We prove that $\Delta B$ is skew-symmetric.
			\begin{align*}
				[\Delta B]_{\textbf i, \textbf j} 
					&= \delta_{\textbf i} b_{\textbf i, \textbf j}\\
					&= d_{\textbf i}|\stab(i)| \frac 1{|\stab(i)|} \sum_{g \in G} a_{g i, j}\\
					&= \sum_{g \in G} d_{gi} a_{gi, j}\\
					&= \sum_{g \in G} d_{j} a_{j,gi}\\
					&= -\sum_{g \in G} d_{gj} a_{gj,i} \\
					&= -|\stab(j)| \frac 1{|\stab(j)|} d_j \sum_{g \in G} a_{g j, i}\\
					&= -\delta_{\textbf j} b_{\textbf j, \textbf i}\\
					&= -[\Delta B]_{\textbf j, \textbf i}.
			\end{align*}
		\end{proof}
		
		\begin{defi}
			Let $Q$ be a valued graph equipped with an admissible group of automorphisms $G$. The \emph{quotient graph} $Q/G$ is the graph associated to the quotient matrix $B_Q/G$. If $Q$ is a quiver, $Q/G$ is called a \emph{quotient quiver}.
		\end{defi}
			
		\begin{exmp}
			Consider again
			$$\xymatrix{Q: 1 \ar@/_10pt/@{--}[rr]& \ar[l] 2 \ar[r] & 3}$$
			with associated matrix
			$$A=\left[\begin{array}{rrr}
				0 & -1 & 0 \\
				1 & 0 & 1\\
				0 & -1 & 0 \\
			\end{array}\right]$$
			and $G=\<(1,3)\>$. We write $\bar 1=\ens{1,3}$ and $\bar 2=\ens 2$ the $G$-orbits. Then,
			$$A/G=\left[\begin{array}{rr}
				0 & -2\\
				1 & 0 \\
			\end{array}\right]$$
			and the associated graph is
			$$\xymatrix{\Delta(A/G): & \bar 1 & \ar[l]_{(2,1)} \bar 2}$$
			which is of type $\B_2$.
		\end{exmp}
	
		We notice that if $A$ is a skew-symmetric matrix such that $\stab_G(i)$ has the same cardinality for every $i \in Q_0$, then $A/G$ is skew-symmetric.
		\begin{exmp}
			We consider the quiver
			$$\xymatrix{
			&&		b_1 \ar@{--}[dd]\\
			A_1:& 	a_1 \ar@{--}[rr] \ar[ru]\ar[rd] && a_2 \ar[ld] \ar[lu]\\
			&&		b_2 
			}$$
			equipped with the admissible automorphism group $G_1=\<(a_1,a_2),(b_1,b_2)\>$. We write $\textbf a=\ens{a_1,a_2}$ and $\textbf b=\ens{b_1,b_2}$, then
			$$\xymatrix{
			A_1/G_1: & \textbf a \ar@<+2pt>[r]\ar@<-2pt>[r] & \textbf b 
			}$$	
			is the Kronecker quiver.
		\end{exmp}
		
		We will also notice that there is no uniqueness in the presentation of a valued graph as a quotient graph.
		\begin{exmp}
			We consider the quiver
			$$
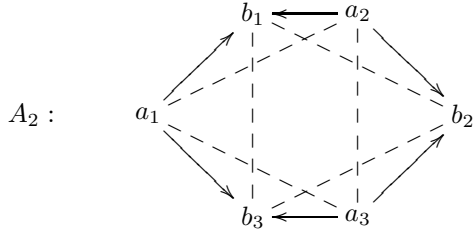
\xymatrix{
			&&		b_1 \ar@{--}[rrd]& \ar[l] a_2 \ar@{--}[dd] \ar[rd] \\
			A_2:& 	a_1 \ar@{--}[rru] \ar[ru]\ar[rd] &&& b_2 \ar@{--}[lld]\\
			&&		b_3 \ar@{--}[uu]& \ar@{--}[llu] a_3 \ar[l] \ar[ru]
			}$$
			equipped with the automorphism group $G_2=\<(a_1,a_2,a_3),(b_1,b_2,b_3)\>$. We set $\textbf a=\ens{a_1,a_2,a_3}$ and $\textbf b=\ens{b_1,b_2,b_3}$, then
			$$\xymatrix{
			A_2/G_2: & \textbf a \ar@<+2pt>[r]\ar@<-2pt>[r] & \textbf b 
			}$$	
			is also the Kronecker quiver.
		\end{exmp}
	\end{subsection}
	
	\begin{subsection}{Cluster algebras and automorphisms}
		We now study the interaction between automorphisms of a skew-symmetrizable matrix $A$ and the cluster algebra associated to the matrix $A$. We fix an admissible pair $(A,G)$ where $A \in M_{Q_0}(\Z)$ is an arbitrary skew-symmetrizable matrix. We set $\textbf u=\ens{u_i \ : \ i \in Q_0}$ and $\textbf v=\ens{v_{\textbf i} \ : \ \textbf i \in \barQ_0}$ to be sets of indeterminates over $\Q$. We set
		$$\mathcal A(A)=\mathcal A(A,\textbf u),$$
		$$\mathcal A(A/G)=\mathcal A(A/G,\textbf v).$$
		
		We define a morphism $\pi$ of $\Z$-algebras, called \emph{projection}:
		$$\pi:\left\{\begin{array}{rcl}
			\Z[\textbf u^{\pm 1}] & \fl & \Z[\textbf v^{\pm 1}]\\
			u_i & \mapsto & v_{\textbf i}.
		\end{array}\right.$$
		
		We define a $G$-action on $\Z[\textbf u^{\pm 1}]$ by setting
		$$g.u_i=u_{gi}$$
		for every $g \in G$, $i \in Q_0$. We now prove that the action of $G$ on $\Z[\textbf u^{\pm 1}]$ induces an action on the cluster algebra $\mathcal A(A)$.
		
		\begin{lem}\label{lem:actsonalgebra}
			Let $A$ be a skew-symmetrizable matrix equipped with a group $G$ of admissible automorphisms. Then $G$ acts on the cluster algebra $\mathcal A(A,\textbf u)$.
		\end{lem}
		\begin{proof}
			We fix a finite set $Q_0$ such that $A \in M_{Q_0}(\Z)$. By definition $g$ acts on $M_{Q_0}(\Z)$ by $g.M=(m_{g^{-1}i, g^{-1}j})$ if $M=(m_{ij})$. 
			For any seed $(B,\textbf v)$ in $\mathcal A(A)$ with $\textbf v=(v_i)_{i \in Q_0}$, we set $g.(B,\textbf v)=(g.B, (v_{g.i})_{i \in Q_0})$.
			
			We prove by induction on $n$ that 
		 	$$g.(\mu_{i_n} \circ \mu_{i_1}(A,\textbf u))=\mu_{g.i_n} \circ \mu_{g.i_1}(A,\textbf u).$$
		 	
		 	We first assume that $n=1$. We fix a vertex $k \in Q_0$. We write $(A',\textbf u')=\mu_k(A,\textbf u)$. Thus $A'=(a'_{ij})$ is given by
		 	$$a'_{ij}=\left\{\begin{array}{ll}
				-a_{ij} & \textrm{if }k \in \ens{i,j}\\
				a_{ij}+\frac{1}{2} \left(|a_{ik}|a_{kj}+a_{ik}|a_{kj}|\right) & \textrm{otherwise}
			\end{array}\right.$$
			and 
			$$u_ku_k'=\prod_{a_ik>0}u_i^{a_{ik}}+\prod_{a_ik<0}u_i^{-a_{ik}}.$$
			
			Thus, $g.(A',\textbf u')=(B,\textbf v)$ is given by $B=(b_j)$ and $\textbf v=(g.u_i')_{i \in Q_0}$ where
			$$b_{ij}=\left\{\begin{array}{ll}
				-a_{g^{-1}i,g^{-1}j} & \textrm{if }k \in \ens{i,j}\\
				a_{g^{-1}i,g^{-1}j}+\frac{1}{2} \left(|a_{g^{-1}i,k}|a_{k,g^{-1}j}+a_{g^{-1}i,k}|a_{k,g^{-1}j}|\right) & \textrm{otherwise}
			\end{array}\right.$$
			so that
			$$b_{ij}=\left\{\begin{array}{ll}
				-a_{i,j} & \textrm{if }gk \in \ens{i,j}\\
				a_{i,j}+\frac{1}{2} \left(|a_{i,gk}|a_{gk,j}+a_{i,gk}|a_{gk,j}|\right) & \textrm{otherwise}
			\end{array}\right.$$
			Also, 
			\begin{align*}
				v_k 	&=gu_k'\\
					&=\frac{1}{g.u_k}\prod_{a_{ik}>0}g.u_i^{a_{ik}} + \prod_{a_{ik}<0}g.u_i^{-a_{ik}}\\
					&=\frac{1}{u_{gk}}\prod_{a_{ik}>0}u_{gi}^{a_{ik}} + \prod_{a_{ik}<0}u_{gi}^{-a_{ik}}\\ \\
					&=\frac{1}{u_{gk}}\prod_{a_{g^{-1}i,g^{-1}k}>0}u_{gi}^{a_{g^{-1}i,g^{-1}k}} + \prod_{a_{g^{-1}i,g^{-1}k}<0}u_{gi}^{-a_{g^{-1}i,g^{-1}k}}\\
			\end{align*}
			and $v_i=g.u_i=u_{gi}$ for $i \neq k$. So that
			$$(B,\textbf v)=\mu_{gk}(A,\textbf ).$$
			
			Assume now that $n>1$. Set
			$$(B,\textbf v)=\mu_{i_{n-1}} \circ \mu_{i_1}(A,\textbf u),$$
			$$(\tilde B,\tilde {\textbf v})=\mu_{g.i_{n-1}} \circ \mu_{g.i_1}(A,\textbf u).$$
			By induction, we now that $g.(B,\textbf v)=(\tilde B,\tilde{\textbf v})$. Thus, we only need to prove that $g.(\mu_{i_n}(B,\textbf v))=\mu_{g.i_n}(\tilde B,\tilde{\textbf v})$. This is proved as in the previous case.
		\end{proof}
	
		\begin{defi}
			A seed $(S, \textbf x)$ in $\mathcal A(A)$ is called \emph{$G$-invariant} if:
			\begin{enumerate}
				\item $gx_i=x_{gi}$ for every $i \in Q_0$, $g \in G$,
				\item $(S,G)$ is an admissible pair.
			\end{enumerate}
		\end{defi}
		
		We will now study the link between the cluster algebra $\mathcal A(A/G)$ and the algebra $\pi(\mathcal A(A))$. For this, we will study a certain class of mutations called \emph{orbit mutations}.
		
		\begin{subsubsection}{Orbit mutations of matrices}
			\begin{lem}\label{lem:coeffAp}
				Let $(A,G)$ be an admissible pair with $A \in M_{Q_0}(\Z)$. Let $\Omega=\ens{i_1, \ldots, i_n}$ be a $G$-orbit in $Q_0$. Then, the coefficients $a_{ij}^{(n)}$ of $A^{(n)}=\mu_{i_n} \circ \cdots \circ \mu_{i_1}(A)$ are given by
				$$a_{ij}^{(n)}=\left\{\begin{array}{ll}
					\displaystyle -a_{ij} & \si i \ou j \in \Omega\\
					\displaystyle a_{ij} + \frac 12  \sum_{k=1}^n \left( |a_{i,i_k}|a_{i_k,j}+a_{i,i_k}|a_{i_k,j}|\right) & \sinon
				\end{array}\right.$$
			\end{lem}
			\begin{proof}
				Set $A=A^{(0)}$ and $A^{(p)}=\mu_{i_p} \circ \cdots \circ \mu_{i_1}(A)$ for every $1 \leq p \leq n$. Mutating at $i_p$, we get
				$$a_{ij}^{(p)}=\left\{\begin{array}{ll}
					\displaystyle -a_{ij}^{(p-1)} & \textrm{if  $i=i_p$ or $j=i_p$} \\
					\displaystyle a_{ij}^{(p-1)} + \frac 12  \left( |a^{(p-1)}_{i,i_p}|a^{(p-1)}_{i_p,j}+a^{(p-1)}_{i,i_p}|a^{(p-1)}_{i_p,j}|\right) & \textrm{otherwise}
				\end{array}\right.$$
	
				By induction, if $i=i_k$ or $j=i_k$ for some $k \leq p-1$, then $a_{ij}^{(p-1)}=-a_{ij}$, thus $a_{ij}^{(p)}=-a_{ij}$ if $i=i_k$ or $j=i_k$ for some $k \leq p$.
	
				Otherwise, 
				\begin{align*}
					a_{ij}^{(p)}
					&=a_{ij}^{(p-1)} + \frac 12  \left( |a^{(p-1)}_{i,i_p}|a^{(p-1)}_{i_p,j}+a^{(p-1)}_{i,i_p}|a^{(p-1)}_{i_p,j}|\right)\\
					&= a_{ij} + \frac 12  \sum_{k=1}^{p-1} \left( |a_{i,i_k}|a_{i_k,j}+a_{i,i_k}|a_{i_k,j}|\right) +\frac 12  \left( |a^{(p-1)}_{i,i_p}|a^{(p-1)}_{i_p,j}+a^{(p-1)}_{i,i_p}|a^{(p-1)}_{i_p,j}|\right)
				\end{align*}
				
				As $a_{i_k,i_l}=0$ for every $k,l$, we get
				$$a^{(p-1)}_{i,i_p}=a_{i,i_p} + \frac 12  \sum_{k=1}^p \left( |a_{i,i_k}|a_{i_k,i_p}+a_{i,i_k}|a_{i_k,i_p}|\right)=a_{i,i_p}$$
				and
				$$a^{(p-1)}_{i_p,j}=a_{i_p,j} + \frac 12  \sum_{k=1}^p \left( |a_{i_p,i_k}|a_{i_k,j}+a_{i_p,i_k}|a_{i_k,j}|\right)=a_{i_p,j}$$
				then
				$$a_{ij}^{(p)}=a_{ij} + \frac 12  \sum_{k=1}^p \left( |a_{i,i_k}|a_{i_k,j}+a_{i,i_k}|a_{i_k,j}|\right).$$
			\end{proof}
			
			\begin{corol}
				Let $(A,G)$ be an admissible pair. Let $\textbf i=\ens{i_1, \ldots, i_n}$ be a $G$-orbit in $Q_0$. Then, for every $\sigma \in \mathfrak S_n$, we have
				$$\mu_{i_n} \circ \cdots \circ \mu_{i_1}(A)=\mu_{i_\sigma(n)} \circ \cdots \circ \mu_{i_\sigma(1)}(A).$$
			\end{corol}
			
			We can thus set the following definition:
			\begin{defi}
				Let $(A,G)$ be an admissible pair, $\textbf i \in \barQ_0$ be a $G$-orbit. Then
				$$\mu_{\textbf i}^G(A):=\left(\prod_{j \in \textbf i} \mu_j\right)(A)$$
				is called the \emph{orbit mutation of $A$ in the direction $\textbf i$}.
			\end{defi}
			
		\end{subsubsection}
		
		\begin{subsubsection}{Orbit mutations of seeds}
			\begin{lem}\label{lem:xin}
				Let $B \in M_{Q_0}(\Z)$ and $G \subset \mathfrak S_{Q_0}$ be a group. Let $(A,\textbf x)$ be a $G$-invariant seed in $\mathcal A(B)$ with $\textbf x=(x_1, \ldots, x_q)$. Let $\Omega=\ens{i_1, \ldots, i_n}$ be a $G$-orbit and $(A^{(n)}, \textbf x^{(n)})=\mu_{i_n}\circ \cdots \circ \mu_{i_1}(A, \textbf x)$. Then, for every $i \in Q_0$, we have
				$$x_i^{(n)}=\left\{\begin{array}{ll}
					\displaystyle x_i 
						& \si i \not \in \Omega\\
					\displaystyle \frac{\prod_{a_{i,i_k}>0} x_i^{a_{i,i_k}}+\prod_{a_{i,i_k}<0} x_i^{-a_{i,i_k}}}{x_{i_k}}  
						& \si i=i_k \textrm{ for some }k\\
				\end{array}\right.$$
			\end{lem}
			\begin{proof}
				For every $p=1, \ldots, n$, we write $(\textbf x^{(p)},A^{(p)})=\mu_{i_p} \circ \cdots \circ \mu_{i_1}(\textbf x,A)$ and $(\textbf x^{(p)},A^{(p)})=(\textbf x,A)$. 
				If $i \not \in \Omega$, for every $p>0$, we have $x_i^{(p)}=x_i^{(p-1)}$ and thus $x_i^{(p)}=x_i$.
				If $i=i_k$ for some $k>0$, then $x_i^{(j)}=x_i$ for every $j<k$ and $x_i^{(l)}=x_i^{(k)}$ for every $l>k$. Exchange relations thus give
				$$x_i^{(k)}=\frac{\prod_{a_{j,i_k}^{(k-1)}>0} (x_j^{(k-1)})^{a_{j,i_k}^{(k-1)}}+\prod_{a_{j,ik}^{(k-1)}<0} (x_j^{(k-1)})^{-a_{j,i_k}^{(k-1)}} }{x_i^{(k-1)}}$$
				But according to lemma  \ref{lem:coeffAp}, $x_i^{(k-1)}=x_i$,  $a_{j,i_k}^{(k-1)}=a_{j,i_k}$. Moreover, $a_{j,i_k}=0$ if $j \in \Omega$, thus if $a_{j,i_k} \neq 0$, we have $x_j^{(k-1)}=x_j$. It follows that
				$$x_i^{(n)}=x_i^{(k)}=\frac{\prod_{a_{j,i_k}>0} x_j^{a_{j,i_k}}+\prod_{a_{j,ik}<0} x_j^{-a_{j,i_k}} }{x_i}$$
				which proves the lemma.
			\end{proof}
			In particular, we see that $(A^{(n)}, \textbf x^{(n)})$ does not depend on the order of the mutations in $\Omega$. We thus set the following definition:
			\begin{defi}
				Let $(A,G)$ be an admissible pair with $A \in M_{Q_0}(\Z)$ and $(S, \textbf x)$ be a $G$-invariant seed in $\mathcal A(A)$. For every $G$-orbit $\textbf i$ in $\barQ_0$, we set
				$$\mu_{\textbf i}^G(S, \textbf x)=\left(\prod_{j \in \textbf i} \mu_j\right)(S, \textbf x)$$
				the \emph{orbit mutation of the seed $(S, \textbf x)$ in direction $\textbf i$}.
			\end{defi}
					
			\begin{corol}\label{corol:orbitmutGinvfaible}
				Let $(B,G)$ be an admissible pair, $(\textbf x,A)$ be a $G$-invariant seed in $\mathcal A(B)$ and $\Omega \in \barQ_0$. Let $(A^{(n)},\textbf x^{(n)})=\mu_\Omega^G((A,\textbf x))$. Then $G$ is an automorphism group of $A^{(n)}$. Moreover, if we denote by $a_{ij}^{(n)}$ the coefficients of $A^{(n)}$ and $\textbf x^{(n)}=(x_i^{(n)}, i \in Q_0)$, we have for all $i \in Q_0$ and $g \in G$:
				$$a_{i,gi}^{(n)}=0 \et x_{gi}^{(n)}=gx_i^{(n)}.$$
			\end{corol}
			\begin{proof}
				We assume that $A \in M_{Q_0}(\Z)$ and $\textbf x=\ens{x_i \ : \ i \in Q_0}$. We write $\Omega=\ens{i_1, \ldots, i_n}$ and we keep notations of lemma \ref{lem:coeffAp}. Let $i,j \in Q_0$. Then, 
				$$g.a^{(n)}_{ij}=a^{(n)}_{g^{-1}i,g^{-1}j}=
				\left\{\begin{array}{ll}
					\displaystyle -a_{g^{-1}i,g^{-1}j} & \si g^{-1}i \textrm{ or } g^{-1}j\in \Omega\\
					\displaystyle \sum_{k=1}^n \frac 12 \left(|a_{g^{-1}i,i_k}|a_{i_k,g^{-1}j}+a_{g^{-1}i,i_k}|a_{i_k,g^{-1}j}| \right)
					& \sinon
				\end{array}\right.$$
				But $g^{-1}i \in \Omega$ \ssi $i \in \Omega$ and in this case $a_{ij}^{(n)}=-a_{ij}=-a_{g^{-1}i,g^{-1}j}=a_{g^{-1}i,g^{-1}j}^{(n)}$.
				But if $i,j \not \in \Omega$, then
				\begin{align*}
				a_{g^{-1}i,g^{-1}j}^{(n)}	
					& =\sum_{k=1}^n \frac 12 \left(|a_{g^{-1}i,i_k}|a_{i_k,g^{-1}j}+a_{g^{-1}i,i_k}|a_{i_k,g^{-1}j}| \right)\\
					&=\sum_{k=1}^n \frac 12 \left(|a_{i,gi_k}|a_{gi_k,j}+a_{i,gi_k}|a_{gi_k,j}| \right)\\
					&=\sum_{k=1}^n \frac 12 \left(|a_{i,i_k}|a_{i_k,j}+a_{i,i_k}|a_{i_k,j}| \right)\\
					&=a_{ij}^{(n)}\\
				\end{align*}
				and thus $G$ is an automorphism group for $A^{(n)}$. Let's prove that $G$ is admissible for $A^{(n)}$. Let $i,j$ be in a same $G$-orbit. We write $i=gj$ for some $g \in G$. According to lemma \ref{lem:coeffAp}, we have
				$$a^{(n)}_{i,j}=a_{ij}+\frac 12 \sum_{k=1}^n \left(|a_{i,i_k}|a_{i_k,j}+a_{i,i_k}|a_{i_k,j}|\right).$$
				Since $(A,G)$ is an admissible pair, we have $a_{i,j}=0$. On the other hand, for every $k=1, \ldots, n$, if $a_{i,i_k}>0$, we necessarily have $a_{i_k,j}\leq 0$. Indeed, if not, there is a path of length 2 from $i$ to $j$ in the graph associated to $A$, which contradicts the fact that $(A,G)$ is an admissible pair. It follows that $a^{(n)}_{i,j}=0$.
	
				Let's now prove that $gx_i^{(n)}=x_{gi}^{(n)}$ for every $i \in Q_0$. If $i \not \in \Omega$, $x_i^{(n)}=x_i$. As $gi \not \in \Omega$, $x_{gi}^{(n)}=x_{gi}$. By hypothesis, $(\textbf x,A)$ is $G$-invariant so that $gx_i=x_{gi}$ and $gx_i^{(n)}=x_{gi}^{(n)}$.
				Now, if $i=i_k$ for some $k \in \ens{1, \ldots, n}$, then
				\begin{align*}
					g.x_i^{(n)}
					&=\frac{\prod_{a_{j,i_k}>0} g.x_j^{a_{j,i_k}}+\prod_{a_{j,ik}<0} g.x_j^{-a_{j,i_k}} }{g.x_i}\\
					&=\frac{\prod_{a_{j,i_k}>0} x_{gj}^{a_{j,i_k}}+\prod_{a_{j,ik}<0} x_{gj}^{-a_{j,i_k}} }{x_{gi}}\\
					&=\frac{\prod_{a_{j,i_k}>0} x_{gj}^{a_{j,i_k}}+\prod_{a_{j,ik}<0} x_{gj}^{-a_{j,i_k}} }{x_{gi}}\\
					&=\frac{\prod_{a_{g^{-1}j,i_k}>0} x_j^{a_{g^{-1}j,i_k}}+\prod_{a_{g^{-1}j,ik}<0} x_j^{-a_{g^{-1}j,i_k}} }{x_{gi}}\\
					&=\frac{\prod_{a_{j,gi_k}>0} x_j^{a_{j,gi_k}}+\prod_{a_{j,gik}<0} x_j^{-a_{j,gi_k}} }{x_{gi}}\\
					&= x_{gi}^{(n)}
				\end{align*}
			\end{proof}
			
			\begin{rmq}\label{rmq:stabilite}
				Beware that in general, if $(A,G)$ is an admissible pair with $A \in M_{Q_0}(\Z)$ and $\Omega \in \overline{Q_0}$, there is no reason for $(\mu^G_{\Omega}(A),G)$ to be a an admissible pair. Indeed, consider the quiver
				$$
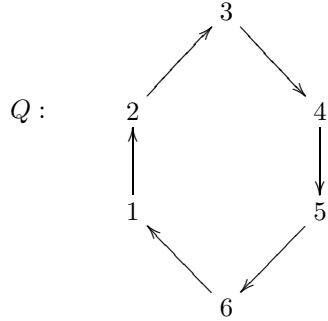
\xymatrix{
						& &3 \ar[rd]\\
					Q: 	& 2 \ar[ru] && 4 \ar[d] \\
						& 1 \ar[u] && 5 \ar[ld] \\
						& &6 \ar[lu]
				}$$
				equipped with the admissible automorphism group $G=\<(1,4)(2,5)(3,6)\>$. Then $\Omega=\ens{3,6}$ is a $G$-orbit and
				$$
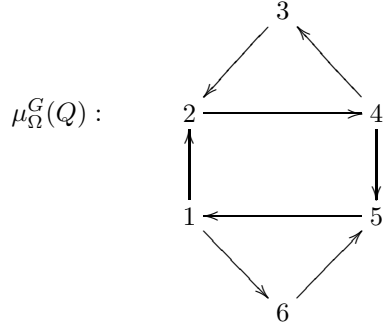
\xymatrix{
							& &3 \ar[ld]\\
				\mu_{\Omega}^G(Q): 	& 2 \ar[rr] && 4 \ar[d]\ar[lu] \\
							& 1 \ar[u]\ar[rd] && 5 \ar[ll] \\
							& &6 \ar[ru]
				}$$
				contains a path of length 2 between from vertex 1 to vertex 4, which both belong to the same $G$-orbit.
			\end{rmq}
			
			We thus introduce the notion of stable admissible pair which will be essential in the following.
			\begin{defi}
				An admissible pair (A,G) is called \emph{stable} if for any finite sequence of $G$-orbits $i_1, \ldots, i_n$,  each of the pairs $(\mu^G_{i_1}(A),G), (\mu^G_{i_2} \circ \mu^G_{i_1}(A),G), \ldots, (\mu^G_{i_n} \circ \cdots \circ \mu^G_{i_1}(A),G)$ is admissible.
			\end{defi}
			Thus, every finite sequence of orbit mutations is well-defined for a stable admissible pair. 
			If $(A,G)$ is a stable admissible pair, we denote by $\Mut^G(A)$ the \emph{orbit mutation class of $A$}, that is, the set of matrices $B$ such that there exists a finite sequence  $(\textbf i_1, \ldots, \textbf i_n)$ of $G$-orbits in $Q_0$ such that
			$$B=\mu_{\textbf i_n}^G \circ \cdots \circ \mu_{\textbf i_1}^G(A).$$
			Note that in particular, we have $\Mut^G(A) \subset \Mut(A)$.
			
			\begin{lem}\label{lem:uneorbitestable}
				Let $(A,G)$ be an admissible pair with $A \in M_{Q_0}(\Z)$. Assume that $G$ has at most one non-trivial orbit in $Q_0$. Then $(A,G)$ is a stable admissible pair.
			\end{lem}
			\begin{proof}
				By induction, it suffices to prove that $(B,G)$ is admissible where $B=\mu^G_{\textbf i}(A)$ for some $G$-orbit $\textbf i$. We denote by $b_{ij}$ the coefficients of $B$. Let $i,j$ be in a same $G$-orbit. According to \ref{corol:orbitmutGinvfaible}, it suffices to prove that for every $k$, $b_{i,k}$ and $b_{j,k}$ are of the same sign. If $k$ is in the $G$-orbit of $i$, then $b_{i,k}=0=b_{j,k}$. Otherwise, $k$ is fixed under the $G$-action by hypothesis. It follows from corollary \ref{corol:orbitmutGinvfaible} that $G$ is an automorphism group of $B$ and thus $b_{i,k}=b_{gi,k}$ for every $g \in G$. In particular, $b_{i,k}$ and $b_{j,k}$ are of the same sign and $(B,G)$ is an admissible pair. Thus, $(A,G)$ is a stable admissible pair.
			\end{proof}
	
			\begin{lem}\label{lem:2orbitesstable}
				Let $(A,G)$ be an admissible pair with $A \in M_{Q_0}(\Z)$. Assume that $G$ has exactly two orbits in $Q_0$. Then, $(A,G)$ is a stable admissible pair.
			\end{lem}
			\begin{proof}
				We denote by $\textbf v=\ens{v_1, \ldots, v_b}$ and $\textbf w=\ens{w_1, \ldots, w_c}$ the two $G$-orbits in $Q_0$ and set $A=(a_{ij})$. Since $(A,G)$ is admissible, we have $a_{v_i,v_j}=0$ for every $i,j \in \ens{1, \ldots, b}$ and $a_{w_i,w_j}=0$ for every $i,j \in \ens{1, \ldots, c}$. Moreover, for every $v \in \textbf v$ and $w \in \textbf w$, $a_{v,w}$ has constant sign. It follows thus from lemma \ref{lem:coeffAp} that $\mu^G_{\textbf v}(A)=-A$ and $\mu^G_{\textbf w}(A)=-A$. In particular, since $(A,G)$ is an admissible pair, $(\mu^G_{\textbf v}(A),G)$ and $(\mu^G_{\textbf w}(A),G)$ are admissible pairs. By induction, $(A,G)$ is thus a stable admissible pair.
			\end{proof}
	
			\begin{prop}\label{prop:Dynkinstable}
				Let $Q$ be a Dynkin quiver equipped with an admissible group of automorphisms $G$. Then, $(Q,G)$ is a stable admissible pair.
			\end{prop}
			\begin{proof}
				Let $(Q,G)$ be an admissible pair with $Q$ Dynkin and let $R \in \Mut^G(Q)$. It follows from corollary \ref{corol:orbitmutGinvfaible} that $R$ is $G$-invariant and that for every $i,j$ in the same $G$-orbit, there is no arrow $i \fl j$ in $R$. It is thus enough to prove that there are no paths of length two between $i$ and $j$ in $R$.
		
				Let $i \in R_0$ and $g \in G$. Assume that there exists a path of length two between $i$ and $gi$. Let's then prove that $R$ contains a minimal cycle of even length. By ``minimal cycle of length $p$'' we mean a cyclic composition of arrows $c=\alpha_1\cdots \alpha_p$ such that any proper sub-path of $c$ is non-cyclic. Assume that there exists a vertex $k \in R_0$ such that $i \fl k \fl gi$. As $G$ is a group of admissible automorphisms of $R$, we know that $k \not \in Gi$. Let
				$$p=\min\ens{n \geq 1 \ : \ g^ni=i}-1,$$
				we thus have a cyclic path
				$$i \fl k \fl gi \fl gk \fl \cdots \fl g^p i \fl g^pk \fl i$$
				in $R$ passing through an even number of vertices.
		
				If all the vertices occurring in this cyclic path are distinct, then we have a minimal cycle of even length. Otherwise, assume that there exists two vertices $v,w$ which are equal among $\ens{g^ni,g^nk, \ : \ 0 \leq n \leq p}$. 
				If $v=g^ni$ and $w=g^mi$ for some $0 \leq m<n \leq p$, we get $i=g^{n-m}i$, which contradicts the minimality of $p$. If $v=g^ni$ and $w=g^mk$ for some $0 \leq m,n \leq p$, we get a contradiction since $k \not \in Gi$. Thus, $v=g^nk$ and $w=g^mk$ for some integers $0 \leq m,n \leq p$. Assume that $n \leq m$, we thus get a new cyclic path:
				$$i \fl k \fl \cdots  \fl g^ni \fl g^nk=g^mk \fl g^{m+1}i \fl \cdots \fl g^p i \fl g^pk \fl i$$
				passing through an even number of vertices. Using inductively this process, we prove that $R$ contains a minimal cycle of even length.
				
				Assume that $Q$ is a Dynkin quiver of type $\A_n$. According to \cite{CCS1}, we know that if $S$ is a quiver in $\Mut(Q)$ then all the minimal cycles in $S$ are 3-cycles. In particular, $R$ does not contain any minimal cycle of length $p$ for every $p \geq 1$ and $(Q,G)$ is thus a stable admissible pair.
				
				Assume that $Q$ is a Dynkin quiver of type $\D_n$ for $n \geq 4$. Then, every automorphism group of $Q$ has at most one non-trivial orbit. It thus follows from lemma \ref{lem:uneorbitestable} that $(Q,G)$ is a stable admissible pair. If $Q$ is of type $\D_3$, then the result follows from lemmas \ref{lem:uneorbitestable} or \ref{lem:2orbitesstable} depending on the considered automorphism group.
				
				Assume that $Q$ is a Dynkin quiver of type $\E_6$. A case by case induction on the mutation class of $Q$ (using for example \cite{Keller:javaapplet}) shows that for any possible automorphism group $G$, $\Mut^G(Q)$ does not contain any quiver such that there is a path of length 2 between two vertices in the same $G$-orbit. In particular, all the pairs $(R,G)$ with $R \in \Mut^G(Q)$ are admissible and thus $(Q,G)$ is stable.
				
				If $Q$ is of type $\E_7$ or $\E_8$, then all the possible automorphism groups of $Q$ are trivial and it thus follows from lemma \ref{lem:uneorbitestable} that $(Q,G)$ is an admissible pair.
			\end{proof}
			
			More generally, Demonet proved that Proposition \ref{prop:Dynkinstable} also holds for any acyclic quiver \cite[Theorem 3.1.11]{Demonet:PhD}:
			\begin{theorem}[\cite{Demonet:PhD}]\label{theorem:acyclicstable}
				Let $Q$ be an acyclic quiver equipped with a group $G$ of admissible automorphisms. Then $(Q,G)$ is a stable admissible pair.
			\end{theorem}
			
		\end{subsubsection}
		
		\begin{subsubsection}{Orbit mutations and mutations of the quotient}
			Let $(A,G)$ be an admissible pair and $(\textbf x, S)$ be a $G$-invariant seed in  $\mathcal A(A)$, we set
			$$\pi(S, \textbf x)=(S/G,\textbf y)$$
			where $\textbf y$ is the $\barQ_0$-tuple given by
			$$y_{\textbf i}=\pi(x_i)$$ 
			for every $\textbf i \in \barQ_0$.
			Then $\pi(S, \textbf x)$ is called \emph{projection of the seed $(S, \textbf x)$}. It is well defined since $(S, \textbf x)$ is $G$-invariant.
	
			\begin{theorem}\label{theorem:mutationquotient}
				Let $(A,G)$ be a stable admissible pair with $A \in M_{Q_0}(\Z)$. Then for every $G$-orbits $\textbf i_1, \ldots, \textbf i_n$ in $\barQ_0$, we have
				$$\mu_{\textbf i_n} \circ \cdots \circ \mu_{\textbf i_1}(\textbf v,A/G) = \pi\left(\mu_{\textbf i_n}^G \circ \cdots \circ \mu_{\textbf i_1}^G(\textbf u,A)\right).$$
			\end{theorem}
			\begin{proof}
				By induction, it suffices to prove it for $n=1$. If $X=(x_{ij})$ is a $G$-invariant matrix, we will denote by $X/G=(\overline{x_{\textbf i, \textbf j}})$ the quotient matrix.
				
				We fix an orbit $\textbf k=\ens{k_1, \ldots, k_n}$ and we write $A^{(n)}=(a_{ij}^{(n)})_{i,j \in Q_0}=\prod_{l \in \textbf k}\mu_l(A)$. According to lemma \ref{lem:coeffAp}, we have
				$$a_{ij}^{(n)}=\left\{\begin{array}{ll}
					\displaystyle -a_{ij} & \si i \ou j \in \textbf k,\\
					\displaystyle a_{ij} + \frac 12  \sum_{s=1}^n \left( |a_{i,k_s}|a_{k_s,j}+a_{i,k_s}|a_{k_s,j}|\right) & \textrm{otherwise.}
				\end{array}\right.$$
				
				On the other hand, we have
				$$\overline{a^{(n)}_{\textbf i, \textbf j}}=\sum_{l \in \textbf i} a_{l,j}^{(n)}$$
				If $\textbf i=\textbf k$ or $\textbf j=\textbf k$, we have
				$$\overline{a^{(n)}_{\textbf i, \textbf j}}=-\sum_{l \in \textbf i} a_{l,j}$$
	
				Otherwise, if $\textbf i, \textbf j \neq \textbf k$, then for every $m \in \textbf i$, 
				$$a_{mj}^{(n)}=a_{mj} + \frac 12  \sum_{s=1}^n \left( |a_{m,k_s}|a_{k_s,j}+a_{m,k_s}|a_{k_s,j}|\right)$$
				and thus
				$$\overline{a^{(n)}_{\textbf i, \textbf j}}=\sum_{m \in \textbf i} \left( a_{mj} + \frac 12  \sum_{s=1}^n \left( |a_{m,k_s}|a_{k_s,j}+a_{m,k_s}|a_{k_s,j}|\right) \right)$$
				or equivalently
				$$\overline{a^{(n)}_{\textbf i, \textbf j}}=\sum_{m \in \textbf i} \left( a_{mj} + \frac 12  \sum_{l \in \textbf k} \left( |a_{m,l}|a_{l,j}+a_{m,l}|a_{l,j}|\right) \right)$$
	
				We now compute the coefficients $(\overline{a_{\textbf i, \textbf j}})'$ of $(\overline A)'=\mu_{\textbf k}(A/G)$.
				By definition of the mutation, 
				$$(\overline{a_{\textbf i, \textbf j}})'=\left\{\begin{array}{ll}
					\displaystyle -\overline a_{\textbf i, \textbf j} & \textrm{if  $\textbf i=\textbf k$ or $\textbf j=\textbf k$,}\\
					\displaystyle \overline a_{\textbf i, \textbf j} + \frac 12  \left(
					|\overline a_{\textbf i, \textbf k}|\overline a_{\textbf k, \textbf j}
					+\overline a_{\textbf i, \textbf k}|\overline a_{\textbf k, \textbf j}|
					\right) & \textrm{otherwise.}
				\end{array}\right.$$
	
				Thus, if $\textbf i=\textbf k$ or $\textbf j=\textbf k$,
				$$(\overline{a_{\textbf i, \textbf j}})'=-\overline{a_{\textbf i, \textbf j}}=-\sum_{l \in \textbf i} a_{l,j}=\overline{a^{(n)}_{\textbf i, \textbf j}}$$
	
				Now, if $\textbf i,\textbf j \neq \textbf k$,
				\begin{align*}
					(\overline{a_{\textbf i,\textbf j}})' 
					&= \overline{a_{\textbf i,\textbf j}}+\frac 12 \left( |\overline{a_{\textbf i,\textbf k}}|\overline{a_{\textbf k,\textbf j}} + \overline{a_{\textbf i,\textbf k}}|\overline{a_{\textbf k,\textbf j}}|\right) \\
					&= \sum_{m \in \textbf i} a_{m,j} + \frac 12 \left( |\sum_{m \in \textbf i}a_{m,k}|\sum_{l \in \textbf k}a_{l,j} + \sum_{m \in \textbf i}a_{m,k}|\sum_{l \in \textbf k}a_{l,j}| \right) \\
					&= \sum_{m \in \textbf i} a_{m,j} + \frac 12 \left( \sum_{m \in \textbf i}|a_{m,k}|\sum_{l \in \textbf k}a_{l,j} + \sum_{m \in \textbf i}a_{m,k}\sum_{l \in \textbf k}|a_{l,j}| \right) \\
					&= \sum_{m \in \textbf i} a_{m,j} + \sum_{m \in \textbf i} \frac 12 \left( |a_{m,k}|\sum_{l \in \textbf k}a_{l,j} + a_{m,k}|\sum_{l \in \textbf k}|a_{l,j}| \right)\\
					&= \sum_{m \in \textbf i} \left[ a_{m,j} + \frac 12 \sum_{l \in \textbf k}\left( |a_{m,k}|a_{l,j} + a_{m,k}||a_{l,j}| \right) \right]\\
				\end{align*}
				and thus
				$$(\overline{a_{\textbf i,\textbf j}})' = \sum_{m \in \textbf i} a_{m,j} + \frac 12 \sum_{m \in \textbf i}\sum_{l \in \textbf k} \left( |a_{m,k}|a_{l,j} + a_{m,k}||a_{l,j}| \right) $$
	
				It thus suffices to prove that
				$$\sum_{m \in \textbf i}\sum_{l\in \textbf k} |a_{m,k}|a_{l,j}+a_{m,k}|a_{l,j}|=\sum_{m \in \textbf i}\sum_{l\in \textbf k} |a_{m,l}|a_{l,j}+a_{m,l}|a_{l,j}|$$
				and thus to prove that
				\begin{equation} \label{doublesum}
				\sum_{m \in \textbf i}\sum_{l\in \textbf k} |a_{m,k}|a_{l,j}=\sum_{m \in \textbf i}\sum_{l\in \textbf k} |a_{m,l}|a_{l,j}
				\end{equation}
				and
				$$\sum_{m \in \textbf i}\sum_{l\in \textbf k} a_{m,k}|a_{l,j}|=\sum_{m \in \textbf i}\sum_{l\in \textbf k} a_{m,l}|a_{l,j}|$$
				In equation (\ref{doublesum}), we consider the coefficient of $a_{lj}$: Let $l_0 \in \textbf k$, in the first sum, the coefficient of $a_{l_0j}$ is $$\sum_{m \in \textbf i} |a_{m,k}|$$ 
				In the second sum, this coefficient is $$\sum_{m \in \textbf i} |a_{m,l_0}|$$
				We write $l_0=g k$ for some $g\in G$. As $a_{m,g k}=a_{g^{-1}m,k}$ and $g^{-1}G=G$, both sums are equal. We identify the same way the coefficients of $|a_{l_0,j}|$ in the second sum and equality (\ref{doublesum}) is proved. Thus, for every $\textbf i, \textbf j \in \barQ_0$, we have
				$$\overline{a_{\textbf i, \textbf j}^{(n)}}=(\overline a_{\textbf i,\textbf j})'$$
				and thus
				$$\mu_{\textbf k}^G(A)/G=\mu_{\textbf k}(A/G).$$
				
				We now get interested in variables. If $j \not \in \ens{k_1, \ldots, k_p}$, $x_j^{(n)}=x_j$ and there is nothing to prove
				According to lemma \ref{lem:xin}, $gx_k^{(n)}=x_{gk}^{(n)}$, thus $\pi(x_{k_j}^{(n)})=\pi(x_{k_s}^{(n)})$ for all $j,s=1, \ldots, n$. 
				Up to reordering, we assume that $k=k_1$ and thus $x_k^{(n)}=x_k^{(1)}=x_k'$. We then write the exchange relation between $x_k$ and $x_k'$ in $\mathcal A(A)$.
				$$x_kx_k'=\prod_{j \in Q_0 \ : \ b_{jk}>0}x_j^{b_{jk}}+\prod_{j \in Q_0 \ : \ b_{jk}<0}x_j^{-b_{jk}}$$
				Applying $\pi$ to this equality, we obtain:
				$$\pi(x_k)\pi(x_k')=\prod_{j \in Q_0 \ : \ b_{jk}>0}\pi(x_j)^{b_{jk}}+\prod_{j \in Q_0 \ : \ b_{jk}<0}\pi(x_j)^{-b_{jk}}$$
				As $(A,G)$ is an admissible pair, $b_{jk}$ and $b_{sk}$ are of the same sign if $j$ and $s$ belong to the same orbit. we can thus group together the terms by signs and we get
				$$\pi(x_k)\pi(x_k')
				=\prod_{\textbf j \in \overline{Q_0} \ : \ b_{jk}>0} \prod_{l \in \textbf j} \pi(x_l)^{b_{lk}}
				+\prod_{\textbf j \in \overline{Q_0} \ : \ b_{jk}<0} \prod_{l \in \textbf j} \pi(x_l)^{-b_{lk}}$$
				As $\pi(x_l)=\pi(x_j)$ for every $l \in \textbf j$, we have
				$$\pi(x_k)\pi(x_k')
				=\prod_{\textbf j \in \overline{Q_0} \ : \ b_{jk}>0} \pi(x_j)^{\sum_{l \in \textbf j} b_{lk}}
				+\prod_{\textbf j \in \overline{Q_0} \ : \ b_{jk}<0} \pi(x_j)^{-\sum_{l \in \textbf j} b_{lk}}$$
				and thus
				$$\pi(x_k)\pi(x_k')
				=\prod_{\textbf j \in \overline{Q_0} \ : \ \overline{b_{\textbf j \textbf k}}>0} \pi(x_j)^{\overline{b_{\textbf j \textbf k}}}
				+\prod_{\textbf j \in \overline{Q_0} \ : \ \overline{b_{\textbf j \textbf k}}<0} \pi(x_j)^{-\overline{b_{\textbf j \textbf k}}}$$
				which is precisely the exchange relation in $\mathcal A(A/G)$ in direction $\textbf k$ with $\pi(x_k)$. If we write $\pi(x_k)'$ the cluster variable in exchange with $\pi(x_k)$ in $\mathcal A(A/G)$, we thus get
				$$\pi(x_{k_j}^{(n)})=\pi(x_k^{(n)})=\pi(x_k')=\pi(x_k)'$$
				for every $j \in \ens{1, \ldots, n}$, which proves the proposition.
			\end{proof}
			
			\begin{corol}\label{corol:corolquotient}
				Let $(A,G)$ be a stable admissible pair. Then every seed in $\mathcal A(A/G)$ is the projection of a $G$-invariant seed in $\mathcal A(A)$. In particular, every cluster variable in $\mathcal A(A/G)$ is the projection of a cluster variable in $\mathcal A(A)$ and thus $\mathcal A(A/G)$ can be identified with a $\Z$-subalgebra of $\pi(\mathcal A(A))$.
			\end{corol}
		\end{subsubsection}
		
		\begin{exmp}\label{exmp:A3toB2}
			We consider again the example of the quiver of type $\A_3$
			$$
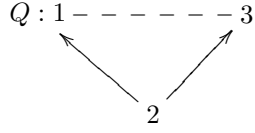
\xymatrix{
			Q: 1 & & 3 \ar@{--}[ll] \\ & \ar[lu] 2 \ar[ru] 
			}$$
			equipped with the admissible group of automorphisms $G=\<(1,3)\>$. Then $(Q,G)$ is a stable admissible pair. The quotient graph of type $\B_2$ is:
			$$\xymatrix{Q/G:&  \overline 2 \ar[r]^{(1,2)} & \overline 1}$$
			where $\overline 1=\ens{1,3}$ et $\overline 2=\ens{2}$. Orbit mutations are thus $\mu_{\overline 1}^G=\mu_1\mu_3$ and $\mu_{\overline 2}^G=\mu_2$. The projection is given by
			$$\pi: \left\{\begin{array}{rcl}
				\Z[u_1^{\pm 1},u_2^{\pm 1},u_3^{\pm 1}] & \fl & \Z[v_{\overline 1}, v_{\overline 2}]\\
				u_1 & \mapsto & v_{\overline 1} \\
				u_2 & \mapsto & v_{\overline 2} \\
				u_3 & \mapsto & v_{\overline 1} \\
			\end{array}\right.$$	
			
			Figure \ref{figure:graphesA3B2} represents the mutation graphs of $\mathcal A(A)$ and $\mathcal A(A/G)$ in the neighbourhood of the initial seeds (where we denoted $v_{\overline i}=v_i$). Vertices with a double circle correspond to seeds obtained by orbit mutations of the initial seed.
			
			Cluster variables in $\mathcal A(Q)$ are given by
			$$\Cl(Q)=\begin{array}{r}
					\left\{\displaystyle u_1,u_2,u_3, 
					\frac{1+u_2}{u_1},\frac{1+u_2}{u_3},
					\frac{1+u_1u_3}{u_2}, \frac{1+u_2+u_1u_3}{u_1u_2},\right.\\
					\left.  \displaystyle \frac{1+u_2+u_1u_3}{u_3u_2},
						\frac{1+2u_2+u_2^2+u_1u_3}{u_1u_2u_3} \right\}
						\end{array}
			$$
			and thus $$\pi(\Cl(Q))=\ens{v_{\overline 1},v_{\overline 2},
					\frac{1+v_{\overline 2}}{v_{\overline 1}}
					\frac{1+v_{\overline 1}^2}{v_{\overline 2}},
					\frac{1+v_{\overline 2}+v_{\overline 1}^2}{v_{\overline 1}v_{\overline 2}},
 					\frac{1+2v_{\overline 2}+v_{\overline 2}^2+v_{\overline 1}^2}{v_{\overline 1}^2v_{\overline 2}}}=\Cl(Q/G)$$
			
			\begin{landscape}
				\begin{figure}
					\setlength{\unitlength}{.45mm}
					\begin{picture}(200,250)(50,100)
						\put(-20,300){$\mathbb A_3$}
							\put(0,300){\circle{5}}
							\put(0,300){\circle{3}}
							\put(0,300){\line(1,0){100}}
							\put(100,300){\circle{5}}
							\put(100,300){\circle{3}}
							\put(100,300){\line(1,1){50}}
							\put(150,350){\circle{5}}
							\put(100,300){\line(1,-1){50}}
							\put(150,250){\circle{5}}
							\put(150,350){\line(0,1){50}}
							\put(200,300){\circle{5}}
							\put(200,300){\circle{3}}
							\put(150,350){\line(1,-1){50}}
							\put(150,250){\line(1,1){50}}
							\put(150,250){\line(0,-1){50}}
						
							\put(200,300){\line(1,0){100}}
							\put(300,300){\circle{5}}
							\put(300,300){\circle{3}}
						\put(-20,290){$\left\{u_1,u_2,u_3\right\}$}
						
						\put(60,290){$\left\{u_1,\frac{1+u_1u_3}{u_2},u_3\right\}$}
						
						\put(110,230){$\left\{\frac{1+u_2+u_1u_3}{u_1u_2},\frac{1+u_1u_3}{u_2},u_3\right\}$}
						
						\put(110,360){$\left\{u_1,\frac{1+u_1u_3}{u_2},\frac{1+u_2+u_1u_3}{u_2u_3}\right\}$}
						
						\put(150,310){$\left\{\frac{1+u_2+u_1u_3}{u_1u_2},\frac{1+u_1u_3}{u_2},\frac{1+u_2+u_1u_3}{u_2u_3}\right\}$}
						
						\put(255,290){$\left\{\frac{1+u_2+u_1u_3}{u_1u_2},\frac{1+2u_2+u_2^2+u_1u_3}{u_1u_2u_3},\frac{1+u_2+u_1u_3}{u_2u_3}\right\}$}

						\put(-20,150){$\mathbb B_2$}
						\put(0,150){\circle{5}}
						\put(0,150){\line(1,0){100}}
						\put(100,150){\circle{5}}
						\put(100,150){\line(1,0){100}}
						\put(200,150){\circle{5}}
						\put(200,150){\line(1,0){100}}
						\put(300,150){\circle{5}}
						
						\put(-20,130){$\left\{v_{\overline 1},v_{\overline 2}\right\}$}
						
						\put(60,130){$\left\{v_{\overline 1},\frac{1+v_{\overline 1}^2}{v_{\overline 2}}\right\}$}
						
						\put(150,130){$\left\{\frac{1+v_{\overline 2}+v_{\overline 1}^2}{v_{\overline 1}v_{\overline 2}},\frac{1+v_{\overline 1}^2}{v_{\overline 2}}\right\}$}
						
						\put(255,130){$\left\{\frac{1+v_{\overline 2}+v_{\overline 1}^2}{v_{\overline 1}v_{\overline 2}},\frac{1+2v_{\overline 2}+v_{\overline 2}^2+v_{\overline 1}v_3}{v_{\overline 1}^2v_{\overline 2}}\right\}$}
					\end{picture}
					\caption{Mutation graphs of $\mathbb A_3$ and $\mathbb B_2$ at the neighbourhood of the initial seeds}\label{figure:graphesA3B2}
				\end{figure}
			\end{landscape}
		\end{exmp}
				
		A consequence of this theorem is that the positivity is preserved. Namely the result is the following:
		\begin{corol}\label{corol:positivitensl}
			Let $(A,G)$ be a stable admissible pair. Assume that every cluster variable in $\mathcal A(A)$ can be written as a Laurent polynomial with positive coefficients in the initial cluster of $\mathcal A(A)$. Then, every cluster variable in $\mathcal A(A/G)$ can be written as a Laurent polynomial with positive coefficients in the initial cluster of $\mathcal A(A/G)$.
		\end{corol}
		\begin{proof}
			We denote by $\textbf v$ the initial cluster of $\mathcal A(A/G)$ and by $\textbf u$ the initial cluster of $\mathcal A(A)$. Let $x$ be a cluster variable in $\mathcal A(Q/G,\textbf v)$. Then $x=\pi(X)$ where $X$ is a cluster variable in $\mathcal A(Q,\textbf u)$ and where $\pi: \Z[\textbf u^{\pm 1}] \fl \Z[\textbf v^{\pm 1}]$ is the morphism of projection, sending $u_i$ on $v_{\textbf i}$ for every $G$-orbit $\textbf i$ of $i$. By assumption, we can write
			$$X=\frac{1}{\textbf u^{\delta(X)}} \sum_\nu a_\nu \textbf u^\nu$$
			where $\nu$ runs over $\N^{Q_0}$ and where $(a_\nu)$ is a family of positive integers with finite support. Thus,
			$$x=\pi(X)=\frac{1}{\pi(u^{\delta(X)})}\sum_\nu a_\nu \pi(\textbf u^\nu)$$
			but for every $\nu \in \N^{Q_0}$, $\pi(\textbf u^\nu)$ is some product $\prod_{\textbf i \in \barQ_0} v_{\textbf i}^{\mu_{\textbf i}}$. In particular, $x$ is a positive Laurent polynomial in $\textbf v$ and the corollary is proved.
		\end{proof}
		
		\begin{exmp}
			The stability condition is essential for theorem \ref{theorem:mutationquotient} to hold. Indeed, consider the quiver
			$$
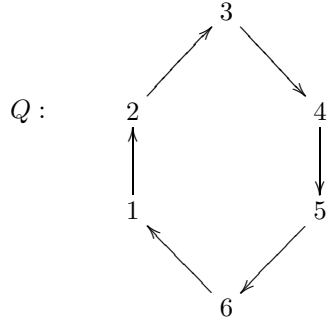
\xymatrix{
						& &3 \ar[rd]\\
					Q: 	& 2 \ar[ru] && 4 \ar[d] \\
						& 1 \ar[u] && 5 \ar[ld] \\
						& &6 \ar[lu]
			}$$
			with the automorphism group $G=\<\prod_{i=1}^3(i;i+3)\> \simeq \Z/2\Z$.	We saw at remark \ref{rmq:stabilite} that $(A,G)$ is not a stable admissible pair. We denote by $\overline i$ the $G$-orbit of a vertex $i$ in $Q_0$. Let $B=B_Q$ be the matrix corresponding to $Q$, $\textbf u=(u_1, \ldots, u_6)$ be the initial seed of $\mathcal A(B)$ and $\textbf v=(v_{\overline 1}, v_{\overline 2}, v_{\overline 3})$ be the initial seed of $\mathcal A(B/G)$. We thus have
			$$\mu^G_{\overline 1}\circ \mu^G_{\overline 2}(\textbf u, B)=(\textbf u',B')$$
			where $$\textbf u'=\left(\frac{u_1u_6+u_3u_6+u_3u_2}{u_1u_2}, \frac{u_1+u_3}{u_2}, u_3, \frac{u_3u_4+u_3u_6+u_6u_5}{u_4u_5}, \frac{u_4+u_6}{u_5}, u_6\right)$$
			and
			$$B'=\left[\begin{array}{rrrrrr}
				0 & 1 & -1 & 0 & 0 & 1\\
				-1 & 0 & 0 & 0 & 0 & 0\\
				1 & 0 & 0 & -1& 0& 0\\
				0& 0& 1& 0& 1& -1\\
				0& 0& 0& -1& 0& 0\\
				-1& 0& 0& 1& 0& 0\\
			\end{array}\right]$$
			Thus,
			$$\pi(\mu^G_{\overline 1}\circ \mu^G_{\overline 2}(\textbf u, B))=
			\left(
				\left(
					\frac{v_{\overline 1}v_{\overline 3}+v_{\overline 3}^2+v_{\overline 2}v_{\overline 3}} {v_{\overline 1}v_{\overline 2}},\frac{v_{\overline 1}+v_{\overline 3}}{v_{\overline 2}},v_{\overline 3}
				\right),
				\left[\begin{array}{rrr}
					0 & 1 & 0 \\
					-1 & 0 & 0 \\
					0 & 0 & 0
				\end{array}\right]
			\right)$$
			but on the other hand,
			$$\mu_{\overline 1}\circ \mu_{\overline 2}(\textbf v, B/G))=
			\left(
				\left(
					\frac{v_{\overline 1}+v_{\overline 3}+v_{\overline 2}}{v_{\overline 1}v_{\overline 2}}, \frac{v_{\overline 1}+v_{\overline 3}}{v_{\overline 2}}, v_{\overline 3}
				\right),
				\left[\begin{array}{rrr}
					0 & 1 & 0 \\
					-1 & 0 & -1 \\
					0 & 1 & 0
				\end{array}\right]
			\right)$$
			and then
			$$\mu_{\overline 1}\circ \mu_{\overline 2}(\textbf v, B/G)) \neq \pi(\mu^G_{\overline 1}\circ \mu^G_{\overline 2}(\textbf u, B)).$$
		\end{exmp}
	\end{subsection}
\end{section}

\begin{section}{$G$-invariant objects in cluster categories}\label{section:interpretationcategorique}
	We will now focus on the case where $(A,G)$ is an admissible pair such that $Q=Q_A$ is an acyclic quiver. When the pair is stable, we get interested in the interpretation we can give of $\mathcal A(A/G)$ in terms of $G$-invariant objects in the cluster category $\CC_Q$ under a certain $G$-action that has to be defined.
	
	\begin{subsection}{$G$-action on the cluster category}
		We first get interested in the representation theory of quivers with automorphisms. A good theoretic framework for this study is the theory of \emph{skew group algebras} developed by Reiten and Riedtmann \cite{RR}. Nevertheless, in our context, we can only consider the more concrete situation of representations of quivers with automorphisms. For general results about these, one can for example refer to the works of Hubery \cite{hubery, huberythesis}.
		
		\begin{subsubsection}{$G$-action on $\rep(Q)$}
			We fix an acyclic quiver $Q$ equipped with a group $G$ of admissible automorphisms. We define an action of $G$ on the category $\rep(Q)$ as follows. Let $V \in \rep(Q)$ and $g \in G$. We define $W=gV$ as the representation given by:
			$$\left\{\begin{array}{lll}
				W(i) &= V(g^{-1}i) & \textrm{ for all } i\in Q_0\\
				W(\alpha) &= V(g^{-1}\alpha) & \textrm{ for all } \alpha\in Q_1\\
			\end{array}\right.$$
			Given two objects $M,N$ in $\rep(Q)$, a morphism of representations $f: M \fl N$ is a family $f=(f_i)_{i \in Q_0}$ such that $f_i \in \Hom_k(M(i),N(i))$ for every $i \in Q_0$. We set
			$$gf=(f_{g^{-1}i})_{i \in Q_0} \in \Hom_{kQ}(gM,gN).$$
			It follows that each $g \in G$ acts functorially on $\rep(Q)$ and that $g^{-1}$ induces a quasi-inverse functor to $g$. In particular, each $g \in G$ defines an auto-equivalence of the category $\rep(Q)$. Moreover, for every $i \in Q_0$, we have $gP_i=P_{gi}$, $gI_i=I_{gi}$ and $gS_i=S_{gi}$.
			
			The action of $G$ on $\rep(Q)$ induces an action of $G$ on the Grothendieck group $K_0(kQ)$. Identifying $K_0(kQ)$ and $\Z^{Q_0}$ with the dimension vector, for every $\textbf d \in \Z^{Q_0}$ and every $g \in G$, We have
			$$g\textbf d=(d_{g^{-1}i})_{i \in Q_0}$$
			and for every representation $M$ in $\rep(Q)$, we thus get
			$$\ddim (gM)=g\ddim M.$$
				
			The following lemma is straightforward:
			\begin{lem}\label{lemgaction}
				Let $M,N \in \rep(Q)$, $\textbf e \in \N^{Q_0}$ and $g \in G$. Then,
				\begin{enumerate}
					\item $\<gM,gN\>=\<M,N\>$,
					\item The variety $\Gr_{g\textbf e}(gM)$ is isomorphic to the variety $\Gr_{\textbf e}(M)$.
				\end{enumerate}
			\end{lem}
		\end{subsubsection}
		
		\begin{subsubsection}{Action on the cluster category}
			The action of $G$ on $\rep(Q)$ induces an action by auto-equivalences on the bounded derived category $D^b(kQ)$ which commutes with $[1]$ and $\tau$. Thus, it induces an action by auto-equivalences on $\CC_Q$. This action is given on the shifts of projective modules by
			$$gP_i[1] \simeq P_{gi}[1].$$
			
			We now prove that the action of $G$ on the objects of the cluster category commutes with the Caldero-Chapoton map. 
			\begin{lem} \label{lem:xgm}
				Let $Q$ be an acyclic quiver equipped with a group $G$ of admissible automorphisms. Then for every object $M$ in $\CC_Q$ and any $g \in G$, we have
				$$X_{gM}=gX_M$$
				
				Equivalently, the following diagram commutes for every $g \in G$:
				$$\xymatrix{
				\Ob(\CC_Q) \ar[r]^{X_?} \ar[d]^g & \Z[\textbf u^{\pm 1}] \ar[d]^g  \\
				\Ob(\CC_Q) \ar[r]^{X_?} & \Z[\textbf u^{\pm 1}]
				}$$
			\end{lem}
			\begin{proof}
				If $M$ decomposes into $M=\bigoplus_i M_i$, we have
				$$X_{gM}=X_{g \left(\bigoplus_i M_i\right)}=X_{\bigoplus_i gM_i}=\prod_i X_{gM_i}.$$
				The action of $G$ on $\Z[\textbf u^{\pm 1}]$ being a morphism of $\Z$-algebras, it suffices to consider the case where $M$ is indecomposable.
		
				As $X_{gP_i[1]}=gu_{i}=u_{gi}=X_{P_{gi}[1]}$, it suffices to prove the result for indecomposable $kQ$-modules. For every modules $M,N$ and every dimension vector $\textbf e \in \N^{Q_0}$, we have
				$$\<gM,gN\>=\<M,N\>$$
				and an isomorphism of varieties $\Gr_{\textbf e}(M) \simeq \Gr_{g{\textbf e}}(gM)$. 
		
				Let $\textbf m= \ddim M$. The $\Z$-linearity of the $G$-action on $\Z[\textbf u^{\pm 1}]$ gives
				\begin{align*}
					gX_M& =\sum_{\textbf e} \chi(\Gr_e(M)) \prod_i u_{gi}^{-\<{\textbf m}, \alpha_i\>-\<\alpha_i, {\textbf m}-{\textbf e}\>} \\
						&= \sum_{\textbf e} \chi(\Gr_e(M)) \prod_i u_i^{-\<{\textbf m}, \alpha_{g^{-1}i}\>-\<\alpha_{g^{-1}i}, {\textbf m}-{\textbf e}\>} 
				\end{align*}
				and
				\begin{align*}
					X_{gM}&=\sum_{\textbf e} \chi(\Gr_e(gM)) \prod_i u_i^{-\<g{\textbf m},\alpha_i\>-\<\alpha_i, g{\textbf m} -{\textbf e}\>} \\
						&= \sum_{\textbf e} \chi(\Gr_{g^{-1}{\textbf e}}(M)) \prod_i u_i^{-\<{\textbf m},\alpha_{g^{-1}i}\>-\<\alpha_{g^{-1}i}, {\textbf m} -g^{-1}{\textbf e}\>} \\
						&= \sum_{\textbf e} \chi(\Gr_{{\textbf e}}(M)) \prod_i u_i^{-\<{\textbf m},\alpha_{g^{-1}i}\>-\<\alpha_{g^{-1}i}, {\textbf m} -{\textbf e}\>} \\
						&= gX_M							
				\end{align*}
			\end{proof}
		\end{subsubsection}
		
		\begin{subsubsection}{$G$-invariant objects in a cluster category}
			If $(Q,G)$ is a stable admissible pair, we saw that the seeds in $\mathcal A(Q/G)$ correspond to some $G$-invariants seeds in $\mathcal A(Q)$. The seeds in $\mathcal A(Q)$ corresponding to rigid objects in $\CC_Q$, it is natural to try to find an interpretation of $G$-invariant seeds in terms of $G$-invariant rigid objects in the cluster category $\CC_Q$.
			
			\begin{defi}
				An object $M$ in $\CC_Q$ is called \emph{$G$-invariant} if for every $g \in G$, we have $gM \simeq M$ in $\CC_Q$. It will be called  \emph{$G$-indecomposable} if it does not have a non-trivial decomposition $M=U \oplus V$ with $U,V$ $G$-invariant objects.
			\end{defi}
	
			\begin{lem}\label{lem:unicitedcp}
				Let $Q$ be an acyclic quiver equipped with a group $G$ of admissible automorphisms. Let $m$, $n$ be two indecomposable objects in $\CC_Q$. Then, $m$ and $n$ are in the same $G$-orbit \ssi 
				$\bigoplus_{U \in Gm}U=\bigoplus_{V \in Gn}V$.
			\end{lem}
			\begin{proof}	
				If $m$ is indecomposable, then $gm$ is indecomposable for every $g \in G$. Thus, the summands appearing in both sums are indecomposable. If we assume that the two sums are isomorphic, as $\CC_Q$ is a Krull-Schmidt category, $n$ appears as an indecomposable summand of both sums and thus $n \in Gm$. The converse is clear.
			\end{proof}
		
			\begin{prop}\label{prop:kqindG}
				Let $M$ be a $G$-indecomposable object in $\CC_Q$. Then, there exists a unique $G$-orbit in $\CC_Q$ containing an object $m$ such that $M=\bigoplus_{U \in Gm}U$. Conversely, for every indecomposable object $m$, the direct sum $\bigoplus_{U \in Gm}U$ is $G$-indecomposable.
			\end{prop}
			\begin{proof}
				Let $m$ be an indecomposable object in $\CC_Q$, we set $M=\bigoplus_{U \in Gm}U$. Then $M$ is $G$-invariant. If $M=N \oplus P$ with $N$ a $G$-indecomposable object and $P \neq 0$, we decompose $N=\bigoplus_i N_i$ into indecomposable objects in $\CC_Q$. The uniqueness of the decomposition of $M$ into indecomposable objects implies that each $N_i$ is in $Gm$. If $N$ is $G$-invariant, then $N$ has all the $gN_i$ as direct summands where $g$ runs over $G$ and thus $N \simeq M$.
				
				Conversely, if $M$ is $G$-indecomposable, we decompose $M=\bigoplus_i M_i$ into indecomposable objects of $\CC_Q$. Each $E_i=\bigoplus_{U \in GM_i}U$ is thus $G$-indecomposable. If $I$ is a set of indices such that $(M_i)_{i\in I}$ is a set of representatives of $G$-orbits of the $M_i$, then $M=\bigoplus_{i \in I} E_i$ is a decomposition into $G$-indecomposable objects. Thus, by hypothesis, $I$ is reduces to one point and all the $M_i$ are in the same orbit. Finally, $M=\bigoplus_{U \in GM_1}U$ where $M_1$ is an indecomposable object in $\CC_Q$. Lemma \ref{lem:unicitedcp} ensures the uniqueness of this decomposition.
			\end{proof}
			
			\begin{corol}
				Every $G$-invariant object $M$ in $\CC_Q$ can be written uniquely (up to reordering and isomorphism) as a direct sum of $G$-indecomposable objects.
			\end{corol}
	
			Thus, we can prove a $G$-invariant analogue to theorem \ref{theorem:correspondanceCK2}.
			\begin{corol}\label{corol:correspondanceGinv}
				Let $Q$ be an acyclic quiver equipped with a group $G$ of admissible automorphisms. Then, $X_?$ induces a 1-1 correspondence from the set of $G$-invariant cluster-tilting objects to the set of $G$-invariant clusters in $\mathcal A(Q)$.
			\end{corol}
			\begin{proof}
				Let $\textbf x=(x_1, \ldots, x_q)$ be a cluster in $\mathcal A(Q)$. According to theorem \ref{theorem:correspondanceCK2}, for every $i=1, \ldots, q$, there exists an unique rigid object $T_i$ in $\CC_Q$ such that $x_i=X_{T_i}$. Moreover, it follows from lemma \ref{lem:xgm} that
				$$\textbf x=g\textbf x=\ens{X_{gT_1}, \ldots, X_{gT_q}}.$$
				The $T_i$ being uniquely determined, it follows that
				$$g \bigoplus_{i=1}^q T_i = \bigoplus_{i=1}^q T_{gi} =\bigoplus_{i=1}^q T_i$$
				and thus $T=\bigoplus_{i=1}^q T_i$ is a $G$-invariant cluster-tilting object in $\CC_Q$.
	
				Conversely, if $T=\bigoplus_{i=1}^q T_i$ is a $G$-invariant cluster-tilting object, then $gT=\bigoplus_{i=1}^q gT_i=T$. It thus follows from lemma \ref{lem:xgm} that
				$\ens{X_{T_1}, \ldots, X_{T_q}}$ is a $G$-invariant cluster.
			\end{proof}
		\end{subsubsection}
	
	\end{subsection}
	
	\begin{subsection}{A denominators theorem}
		We now give a categorical interpretation for the denominators of cluster variables in the cluster algebra $\mathcal A(Q/G)$. If $Q$ is an acyclic quiver with a group $G$ of admissible automorphisms, we define the \emph{projection} on the Grothendieck group $\pi: K_0(kQ) \fl \Z^{\barQ_0}$ by setting
		$$\pi\left( (d_i)_{i \in Q_0} \right)=(\sum_{j \in \textbf i}d_j)_{\textbf  i \in \barQ_0}.$$
		
		\begin{theorem} \label{theorem:theodenGinv}
			Let $Q$ be an acyclic quiver equipped with a group $G$ of admissible automorphisms. Let $m$ be an object in $\CC_Q$. Then,
			$$\delta(\pi(X_m))=\pi(\ddim m).$$
			Equivalently, the following diagram commutes:
			$$\xymatrix{
				\Ob(\mathcal C_Q) \ar[r]^{X_?} \ar[rd]^{\pi \circ \ddim} & \mathcal \Z[u_i^{\pm 1}: i \in Q_0]  \ar[d]^{\delta \circ \pi} \\
				& \Z^{\barQ_0}
			}$$
		\end{theorem}
		
		Before proving the theorem, we prove the following technical lemma:
		\begin{lem}\label{lem:lemtheodenGinv}
			Let $M$ be a $G$-invariant $kQ$-module. Then,
			\begin{enumerate}
				\item for every $l \in Q_0$ and any submodule $N \subset M$, 
					$$(\ddim M)_l \leq <N,S_l>+<S_l, M/N>;$$
				\item for every $l \in Q_0$, there exists a submodule $N \subset M$ such that for every $k \in  \textbf l$, we have $$(\ddim M)_k=<N,S_k>+<S_k, M/N>.$$
			\end{enumerate}
		\end{lem}
		\begin{proof}
			The first point is proved in \cite{CK2}. For the second point, we write $\textbf d =\ddim M$, $\textbf e = \ddim N$ where $N$ is a submodule of $M$. Let $n_k=-<N,S_k>-<S_k, M/N>$. For every $k \in G.l$, we have
			$$n_k=-d_k+\sum_{i \fl k} e_i+\sum_{k \fl j}(d_j - e_j).$$
			
			Let $N$ be the submodule of $M$ generated by the sum of the $M_j$ such that there exists $k \in G.l$ with $k \fl j$ in $Q_1$. It suffices to prove that for any $k \in G.l$, we have $n_k=-d_k$. For this, it is enough to prove that both sums are zero.
			
			Fix $k \in G.l$. Assume that the first sum is non-zero. Then, there exists some $i$ such that $i \fl k$ in $Q_1$ and $e_i \neq 0$. If $e_i \neq 0$, then by definition of $N$, there exists $j \in \overline l$ and a path in $Q$ from $j$ to $i$ and thus $j \fl \cdots \fl i \fl k$ in $Q$. As $j$ and $k$ are in the same $G$-orbit and that each element in $G$ has a finite order, it follows that there exists a cycle in $Q$, which is a contradiction. The second sum is zero by construction. 
		\end{proof}
		
		We now prove theorem \ref{theorem:theodenGinv}:
		\begin{proof}\emph{(of theorem \ref{theorem:theodenGinv})}
			We first notice that for any objects $m,n$ in $\CC_Q$, we have 
			$$\begin{array}{rl}
				\delta(\pi(X_{m \oplus n})) &= \delta(\pi(X_m)\pi(X_n)) \\
					&= \delta(\pi(X_m))+\delta(\pi(X_n))
			\end{array}$$
			we can thus assume that $m$ is indecomposable. If $m=P_i[1]$ for some $i \in Q_0$, the result holds. It is thus enough to consider the case where $m$ is a $kQ$-module.
			
			Let $M=\bigoplus_{U \in Gm} U$, we have
			$$X_M=\prod_{U \in Gm}X_U$$
			According to lemma \ref{lem:xgm}, $X_{gm}=gX_m$ for every $g \in G$ and thus
			$$\pi(X_M)=\pi(X_m)^{|Gm|}$$
			and $\delta(\pi(X_M))=|Gm|\delta(\pi(X_m))$.
			
			The variable $X_M$ is a sum of monomials $\prod x_k^{n_k}$ where $n_k=-<N,S_k>-<S_k, M/N>$ for a certain submodule $N$ of $M$. Fix $i \in Q_0$, lemma \ref{lem:lemtheodenGinv} proves that there exists a term in the expansion of $X_M$ such that the exponent of $u_j$ in the denominator of this monomial is maximal and equal to $\dim M(j)$ for every $j \in \textbf i$. Thus, under projection, the exponent of $v_{\textbf  i}$ in $\pi(X_M)=\pi(X_m)^{|Gm|}$ is $\sum_{j \in \textbf i}\dim M(j)$.
			
			Thus,
			\begin{align*}
				\delta(\pi(X_M)) 
					&= \left( \sum_{j \in \overline i} \dim M(j) \right)_{\overline i \in \barQ_0}\\
					&= \pi(\ddim M) \\
					&= \pi\left(\ddim \bigoplus_{U \in Gm} U \right)\\
					&= \pi \left(\sum_{U \in Gm} \ddim U \right)\\
					&= \sum_{U \in Gm} \pi \left(\ddim U \right)\\
					&= |Gm| \pi (\ddim m)
			\end{align*}		
			It follows that $\delta(\pi(X_m))=\pi (\ddim m)$, which proves the theorem.
		\end{proof}
	\end{subsection}
	
	\begin{subsection}{The finite type case}	
		We now get interested in the case where $Q$ is a Dynkin quiver. We prove in particular that if $Q$ is a Dynkin quiver equipped with a group $G$ of automorphisms, then $\mathcal A(Q/G)=\pi(\mathcal A(Q))$. Moreover, we prove that the denominator vector induces a 1-1 correspondence from the set of cluster variables in $\mathcal A(Q/G)$ to the set of almost positive roots of $Q/G$.
		
		Let $Q$ be a Dynkin quiver and $G$ is a group of automorphisms of $G$. Since $Q$ is acyclic, $G$ is necessarily admissible and moreover it follows from proposition \ref{prop:Dynkinstable} that $(Q,G)$ is a stable admissible pair. We recall in figures \ref{figure:AntoBn}, \ref{figure:DntoCn}, \ref{figure:E6toF4}, \ref{figure:D4toG2} below all the non-trivial quotients $Q/G$ when $Q$ is Dynkin (see also \cite{Lusztig:quantumbook}).
		
		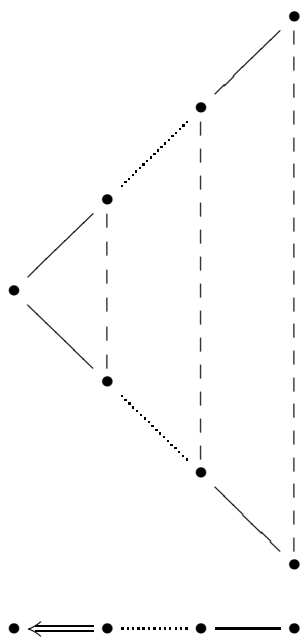
\begin{figure}[H]
			$$\xymatrix{
			&&&\bullet \ar@{--}[dddddd]\\
			&& \bullet \ar@{-}[ru] \ar@{--}[dddd]\\
			& \bullet \ar@{.}[ru] \ar@{--}[dd]\\
			\bullet \ar@{-}[ru] \ar@{-}[rd]\\
			& \bullet \ar@{.}[rd]\\
			&& \bullet \ar@{-}[rd]\\
			&&& \bullet
			}$$ 
			
			$$\xymatrix{
			\bullet & \ar@{=>}[l] \bullet & \bullet \ar@{.}[l] & \bullet \ar@{-}[l] 
			}$$
			\caption{Quotient from $A_{2n-1}$ to $C_n$}\label{figure:AntoBn}
		\end{figure}
	
		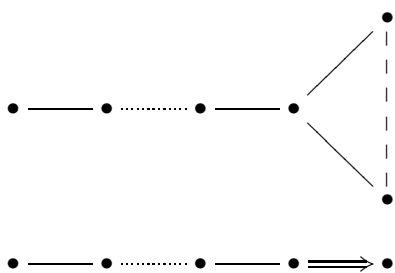
\begin{figure}[H]
			$$\xymatrix{
			&&&& \bullet \ar@{--}[dd]\\
			\bullet \ar@{-}[r] & \bullet \ar@{.}[r] & \bullet \ar@{-}[r] & \bullet \ar@{-}[ru]\ar@{-}[rd]\\
			&&&& \bullet \\
			}$$ 
			
			$$\xymatrix{
			\bullet & \ar@{-}[l] \bullet & \bullet \ar@{.}[l] & \bullet \ar@{-}[l] \ar@{=>}[r] & \bullet
			}$$
			\caption{Quotient from $D_{n+1}$ to $B_n$}\label{figure:DntoCn}
		\end{figure}
		
		\begin{figure}[H]
			$$\xymatrix{
			&&& \bullet \ar@{--}[dddd]\\
			&& \bullet \ar@{-}[ru] \ar@{--}[dd]\\
			\bullet \ar@{-}[r] & \bullet \ar@{-}[ru] \ar@{-}[rd]\\
			&& \bullet \ar@{-}[rd] \\
			&&& \bullet\\
			}$$ 
			
			$$\xymatrix{
			\bullet & \ar@{-}[l]  \bullet & \bullet \ar@{=>}[l] & \bullet \ar@{-}[l]
			}$$
			\caption{Quotient from $E_{6}$ to $F_4$}\label{figure:E6toF4}
		\end{figure}
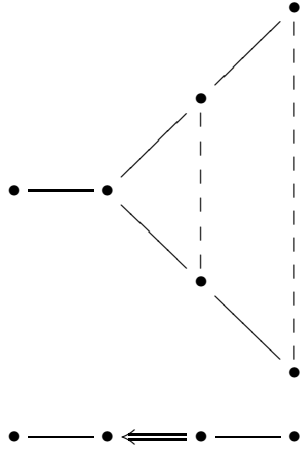
	
		\begin{figure}[H]
			$$\xymatrix{
			&& \bullet \ar@/^/@{--}[dd] \ar@/_2pt/@{--}[lld]\\
			\bullet \ar@{-}[r] & \bullet \ar@{-}[ru] \ar@{-}[rd]\\
			&& \bullet \ar@/^2pt/@{--}[llu]\\
			}$$ 
			
			$$\xymatrix{
			\bullet & \bullet \ar@{-}[l]^{(3,1)}
			}$$
			\caption{Quotient from $D_{4}$ to $G_2$}\label{figure:D4toG2}
		\end{figure}
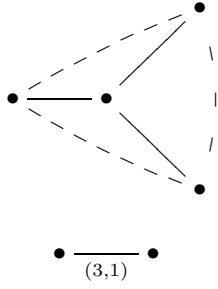
	
		The following lemma can be obtained by direct computation:
		\begin{lem}\label{sysracines}
			Let $Q$ be a quiver of Dynkin type equipped with a group $G$ of automorphisms and let $Q/G$ be the quotient graph. Then,
			$$\pi(\Phi_{\geq -1}(Q))=\Phi_{\geq -1}(Q/G).$$
		\end{lem}
		
		\begin{rmq}\label{rmq:Huberyreflexions}
			According to Hubery's works (see \cite{hubery}), we can prove that $\Phi_{\geq -1}^{\textrm{re}}(Q/G) \subset \pi(\Phi^{\textrm{re}}_{\geq -1}(Q))$ for any acyclic quiver $Q$. Methods used in this context are similar to those we used above. Namely, instead of considering orbit mutations, the author considers orbit reflections in the Weyl group. More precisely, if $\ens{\alpha_i \ : \ i \in Q_0}=\Pi(Q)$ is the set of simple roots of $Q$, for every $G$-orbit $\textbf i$, we set $\beta_{\overline i}=\pi(\alpha_i)=\sum_{j \in \textbf i}\alpha_j$. We denote by $s_{i}$ the reflections in the Weyl group of $Q$ and $s_{\overline i}$ the reflections in the Weyl group of $Q/G$. As $G$ is admissible $Q$, reflections taken in vertices in a same $G$-orbit commute and we define the \emph{orbit reflection} $r_{\textbf i}=\prod_{j \in \textbf i} s_j$. It thus follows from \cite{hubery} that $\pi \circ r_{\overline i}=\sigma_{\overline i} \circ \pi$ for every $G$-orbit $\textbf i$. Thus, by induction from the simple roots, we prove that the real roots of $Q/G$ can be obtained as projections of certain real roots of $Q$.
		\end{rmq}
		
		In order to prove that for a finite type quiver $Q$, we have equality between $\pi(\mathcal A(Q))$ and $\mathcal A(Q/G)$, we will need the following lemma:
		\begin{lem}\label{lem:domainefondamental}
			Let $Q$ be a Dynkin quiver equipped with a group $G$ of automorphisms. Let $\gamma, \beta \in \Phi(Q)$ be roots such that $\pi(\gamma)=\pi(\beta)$. Then there exists $g \in G$ such that $\gamma=g.\beta$.
		\end{lem}
		\begin{proof}
			Let $\Pi(Q)=\ens{\alpha_i, i \in Q_0}$ be the set of simple roots of $Q$ and $\Pi(Q/G)=\ens{\alpha_{\textbf i}, \textbf i \in \textbf Q_0}$ be the set of simple roots of $Q/G$.
			Let $\gamma, \beta$ be roots in $\Phi(Q)$. We write $\gamma=\sum_{i \in Q_0} \gamma_i \alpha_i$ and $\beta=\sum_{i \in Q_0} \beta_i \alpha_i$. By hypothesis, for every $\textbf i \in \textbf Q_0$, we have $\sum_{k \in \textbf i}\beta_k=\sum_{k \in \textbf i}\gamma_k$. If $\gamma$ is a simple root, then $\beta$ is also a simple root and there exists $g \in G$ such that $\gamma=g.\beta$. Otherwise, according to lemma \ref{sysracines} $\pi(\beta)$ is an element in $\Phi(Q/G)$. With notations of remark \ref{rmq:Huberyreflexions}, there exists a sequence of vertices $\textbf i_1, \ldots, \textbf i_n$ in $\textbf Q_0$ such that $\sigma_{\textbf i_n} \circ \cdots \circ \sigma_{\textbf i_1}(\pi(\beta))$ is a simple root in $Q/G$. Then, it follows from remark \ref{rmq:Huberyreflexions} that
			\begin{align*}
				\pi(r_{\textbf i_n} \circ \cdots \circ r_{\textbf i_1}(\beta))
					&= \sigma_{\textbf i_n} \circ \cdots \circ \sigma_{\textbf i_1}(\pi(\beta))\\
					&= \sigma_{\textbf i_n} \circ \cdots \circ \sigma_{\textbf i_1}(\pi(\gamma))\\
					&= \pi(r_{\textbf i_n} \circ \cdots \circ r_{\textbf i_1}(\gamma)).
			\end{align*}
			Then, it follows from the discussion in the simple case that there exists some $g \in G$ such that 
			$$g.r_{\textbf i_n} \circ \cdots \circ r_{\textbf i_1}(\beta)=r_{\textbf i_n} \circ \cdots \circ r_{\textbf i_1}(\gamma).$$
			But it follows from \cite{huberythesis} that the action of $g$ commutes with the $r_{\textbf i}$ and thus
			$$r_{\textbf i_n} \circ \cdots \circ r_{\textbf i_1}(g.\beta)=r_{\textbf i_n} \circ \cdots \circ r_{\textbf i_1}(\gamma)$$
			and $\gamma=g.\beta$.
		\end{proof}
		
		\begin{theorem}\label{theorem:finitequotient}
			Let $Q$ be a Dynkin quiver equipped with a group $G$ of automorphisms and let $Q/G$ be the quotient graph. Then
			$$\pi(\Cl(Q)) = \Cl(Q/G).$$
			In particular,
			$$\pi(\mathcal A(Q))=\mathcal A(Q/G).$$
		\end{theorem}
		\begin{proof}
			Since $Q$ is Dynkin, $(Q,G)$ is a stable admissible pair and thus it follows from corollary \ref{corol:corolquotient} that we have an inclusion $\iota: \Cl(Q/G) \fl \pi(\Cl(Q))$. Since we know (see for example \cite{Zhu:applications}) that denominator vectors induce 1-1 correspondences
			$$\delta_1: \Cl(Q/G) \xrightarrow{\sim} \Phi_{\geq -1}(Q/G),$$
			$$\delta_3: \Cl(Q) \xrightarrow{\sim} \Phi_{\geq -1}(Q).$$
			We have a commutative diagram
			$$\xymatrix{
				\Cl(Q/G) \ar[r]^{\iota} \ar[d]_{\delta_1}^{\sim} & \pi(\Cl(Q)) \ar[d]_{\delta_2} & \ar[l] \Cl(Q)\ar[d]_{\delta_3}^{\sim}\\
				\Phi_{\geq -1}(Q/G) \ar@{=}[r] & \pi(\Phi_{\geq -1}(Q)) & \ar[l] \Phi_{\geq -1}(Q)
			}$$
			where the square on the left hand side clearly commutes and the square on the right hand side commutes according to theorem \ref{theorem:theodenGinv}. Since, $\delta_2 \circ \iota=\delta_1$ is a 1-1 correspondence, $\delta_2$ is surjective. Let's now prove that $\delta_2$ is injective. Fix $x,y \in \Cl(Q)$ such that $\delta_2(\pi(x))=\delta_2(\pi(y))$. It follows from \ref{theorem:theodenGinv} that $\pi(\delta_3(x))=\pi(\delta_3(y))$. According to lemma \ref{lem:domainefondamental}, there exists $g \in G$ such that $\delta_3(y)=g.\delta_3(x)$. As $\delta_3(y)$ is an element in $\Phi_{\geq -1}(Q)$, there exists an unique element $\Cl(Q)$ with denominator vector $\delta_3(y)$ and thus according to lemma \ref{lem:xgm}, this element is $g.x$. Thus, $y=g.x$ and $\pi(y)=\pi(x)$. It follows that $\delta_2: \pi(\Cl(Q)) \fl \pi(\Phi_{\geq -1}(Q))$ is a 1-1 correspondence and thus $\iota$ also. Thus, we have $\pi(\Cl(Q))=\Cl(Q/G)$, which proves the theorem.
		\end{proof}
	\end{subsection}
	
	\begin{subsection}{The affine case}
		We proved that $\pi(\Cl(Q))=\Cl(Q/G)$ if $Q$ is a Dynkin quiver equipped with a group $G$ of automorphisms. We now prove that if $Q$ is of infinite representation type, then the inclusion $\Cl(Q/G) \subset \pi(\Cl(Q))$ can be proper.
			
		Consider the quiver $Q$ of affine type $\Daffine_4$
		$$
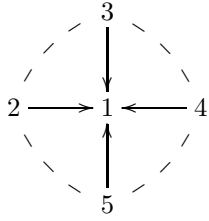
\xymatrix{
			& 3 \ar[d] \ar@{--}@/^/[rd]\ar@{--}@/_/[ld]\\
		2 \ar[r] & 1 & \ar[l] 4\\
			& 5 \ar[u] \ar@{--}@/_/[ru]\ar@{--}@/^/[lu]
		} $$
		equipped with the admissible group of automorphisms $\<g\>$ where $g$ is the 4-cycle $g=(2345) \in \mathfrak S_5$. It follows from lemma \ref{lem:2orbitesstable} that $(Q,G)$ is a stable admissible pair.
		
		The quotient graph $Q/G$ is the valued graph of type $\A_2^{(2)}$ given by
		$$\xymatrix{
		\textbf 1 & \ar[l]_{(4,1)} \textbf 2 
		}$$
	
		$Q$ and $Q/G$ are both of affine types. The minimal positive imaginary roots are respectively $\delta_Q=(21111)$ and $\delta_{Q/G}=(1,2)$.
		
		For every $i \neq j$ in $\ens{2,3,4,5}$, we set $M_{ij}$ to be the unique (up to isomorphism) indecomposable representation of dimension vector $\alpha_1 + \alpha_i + \alpha_j$. Then the $M_{ij}$ are quasi-simple modules in exceptional tubes and $\tau M_{ij}=M_{rs}$ where $\ens{r,s}=\ens{2,3,4,5} \setminus \ens{i,j}$. In particular, we have $\Ext^1(M_{ij},M_{rs}) \neq 0$.
		
		For every $i \neq j \in \ens{2,3,4,5}$, $\ddim M_{ij}$ is a real Schur root and thus $\Ext^1_{kQ}(M_{ij},M_{ij})=0$. Theorem \ref{theorem:correspondanceCK2} implies that $X_{M_{ij}}$ is a cluster variable in $\mathcal A(Q)$.
		
		We now consider the representation
		$$
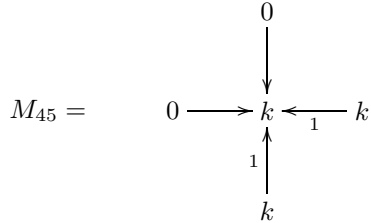
\xymatrix{
			&& 0 \ar[d] \\
		M_{45}= &0 \ar[r] & k & \ar[l]^{1} k\\
			&& k \ar[u]^1 
		}$$
		If $\pi(X_{M_{45}})$ is a cluster variable in $\mathcal A(Q/G)$, then according to corollary \ref{corol:corolquotient}, there exists a cluster $\textbf x=\ens{x_1, \ldots, x_5}$ obtained after a sequence of orbit mutations such that $\pi(X_{M_{45}})=\pi(x_{i_0})$ for some $i_0 \in \ens{1, \ldots, 5}$. According to theorem \ref{theorem:correspondanceCK2}, there exists an indecomposable rigid object $N$ such that $x_{i_0}=X_N$. But according to theorem \ref{theorem:theodenGinv}, $\pi(\ddim N)=\pi(\ddim M_{45})=(1,2)$, so $N$ is necessarily one of the $M_{ij}$. According to lemma \ref{lem:xin}, we know moreover that $\textbf x$ is a $G$-invariant cluster and thus it follows from corollary \ref{corol:correspondanceGinv} that $M_{ij}$ is a direct summand of a $G$-invariant cluster-tilting object in $\CC_Q$. But $\Ext^1_{kQ}(g^2M_{ij}, M_{ij}) \neq 0$, hence a contradiction. Thus, $\pi(X_{M_{45}}) \in \pi(\Cl(Q)) \setminus \Cl(Q/G)$ and
		$$\Cl(Q/G) \subsetneq \pi(\Cl(Q)).$$
	\end{subsection} 
\end{section}

\begin{section}{Mutation-finite diagrams}\label{section:mutationfinitude}
	A skew-symmetrizable matrix $A$ is called \emph{mutation-finite} if $\Mut(A)$ is finite. The problem of classifying all mutation-finite matrices is still open and has been subject to various developments \cite{FST:finitemutationtype, cluster2, FST, Hille:cluster, seven:symmetrizable, BR:extended, Derksen:finite}. In the skew-symmetric case, the problem was solved by Felikson, Shapiro and Tumarkin in \cite{FST:finitemutationtype} using a list proposed by Fomin, Shapiro, Thurston, Derksen and Owen \cite{FST, Derksen:finite}. In the skew-symmetrizable case, the main result is due to Seven. It states that an acyclic valued graph is mutation-finite if and only if it is an affine valued graph, generalizing a theorem of Buan and Reiten \cite{BR:extended} for the symmetric case. In this section, realizing affine valued graphs as quotients of affine quivers, we prove independently of Seven's result that certain affine diagrams are mutation-finite.
	
	In this section, the key result will be theorem of Buan and Reiten:
	\begin{theorem}[\cite{BR:extended}]\label{theorem:mutationfiniteBR}
		Let $Q$ be an acyclic quiver with at least three vertices. Then $Q$ is mutation-finite if and only if $Q$ is of Dynkin or affine type.
	\end{theorem}	
	
	We recall that an \emph{non-simply laced affine graph} is an orientation of on of the diagrams listed in figure \ref{figure:nonsimplyaffine}.
	\begin{figure}[H]
		$$\xymatrix{
		\tilde {\mathbb A}^{(2)}_1:  & \bullet \ar@{-}[r]^{(1,4)} & \bullet
		}$$
		$$\xymatrix{
			\tilde{\mathbb B}_n (n\geq 2): & \bullet \ar@{-}[r]^{(1,2)} & \bullet \ar@{-}[r] & \bullet \ar@{--}[r]& \bullet \ar@{-}[r] & \bullet \ar@{-}[r]^{(2,1)} & \bullet\\
		}$$
		$$\xymatrix{
			\tilde{\mathbb C}_n (n\geq 2): & \bullet \ar@{-}[r]^{(2,1)} & \bullet \ar@{-}[r] & \bullet \ar@{--}[r]& \bullet \ar@{-}[r] & \bullet \ar@{-}[r]^{(1,2)} & \bullet\\
		}$$
		$$\xymatrix{
			\tilde{\mathbb {BC}}_n (n\geq 2): & \bullet \ar@{-}[r]^{(2,1)} & \bullet \ar@{-}[r] & \bullet \ar@{--}[r]& \bullet \ar@{-}[r] & \bullet \ar@{-}[r]^{(1,2)} & \bullet\\
		}$$
		$$\xymatrix{
			&\bullet \ar@{-}[rd]\\
			\tilde{\mathbb {BD}}_n (n\geq 2): &&  \bullet \ar@{-}[r] & \bullet \ar@{--}[r]& \bullet \ar@{-}[r] & \bullet \ar@{-}[r]^{(2,1)} & \bullet\\
			&\bullet \ar@{-}[ru]
		}$$
		$$\xymatrix{
			&\bullet \ar@{-}[rd]\\
			\tilde{\mathbb {CD}}_n (n\geq 2): &&  \bullet \ar@{-}[r] & \bullet \ar@{--}[r]& \bullet \ar@{-}[r] & \bullet \ar@{-}[r]^{(1,2)} & \bullet\\
			&\bullet \ar@{-}[ru]
		}$$
		$$\xymatrix{\tilde{\mathbb F}^{(1)}_{4}:  & \bullet \ar@{-}[r] & \bullet \ar@{-}[r] & \bullet \ar@{-}[r]^{(1,2)} &\bullet \ar@{-}[r] &\bullet }$$
		$$\xymatrix{\tilde{\mathbb F}^{(2)}_{4}:  & \bullet \ar@{-}[r] & \bullet \ar@{-}[r] & \bullet \ar@{-}[r]^{(2,1)} &\bullet \ar@{-}[r] &\bullet }$$
		$$\xymatrix{\tilde{\mathbb G}^{(1)}_{2}:  & \bullet \ar@{-}[r] & \bullet \ar@{-}[r]^{(1,3)} &\bullet}$$
		$$\xymatrix{\tilde{\mathbb G}^{(2)}_{2}:  & \bullet \ar@{-}[r] & \bullet \ar@{-}[r]^{(3,1)} &\bullet}$$
	\caption{non-simply laced diagrams}\label{figure:nonsimplyaffine}
	\end{figure}
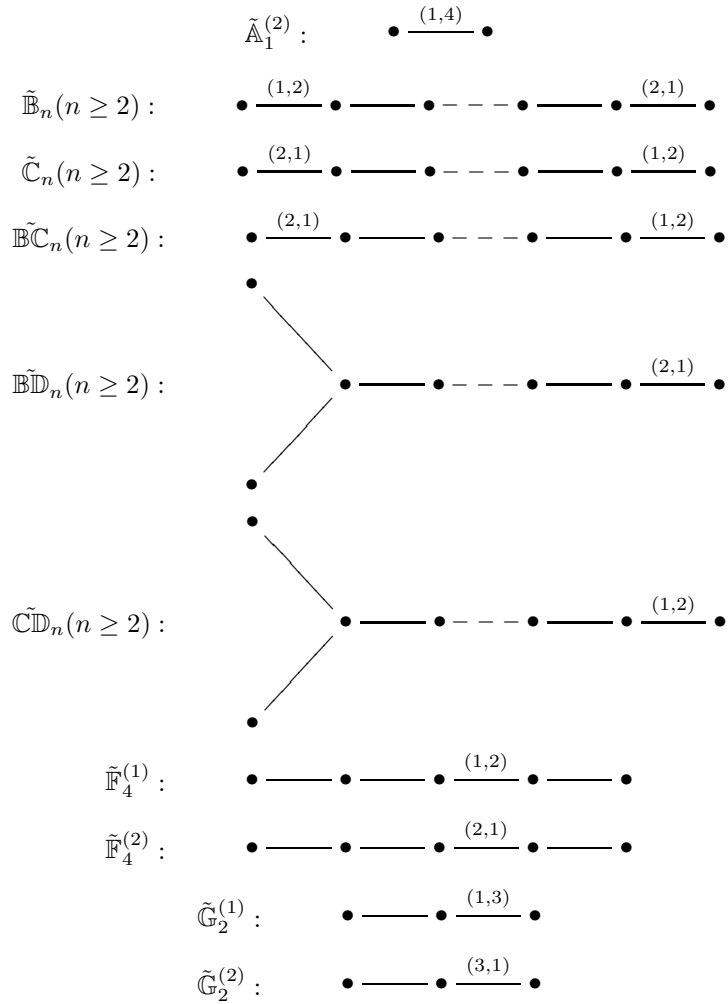

	Lusztig proved that any affine non-simply-laced diagram can be realized as a quotient of a simply-laced affine by a group of automorphisms. More precisely,
	\begin{lem}[\cite{Lusztig:quantumbook}]\label{lem:affinequotient}
		Let $\Delta$ be a non-simply-laced valued graph. Then there exists an affine quiver $Q$ equipped with a group $G$ of automorphisms such that $\Delta=Q/G$.
	\end{lem}
	The list of possible quotients in the Dynkin case is well-known. In the case of affine types, we recall the  possible quotients in figures \ref{figure:D4atoA21a}-\ref{figure:E6atoG22}.

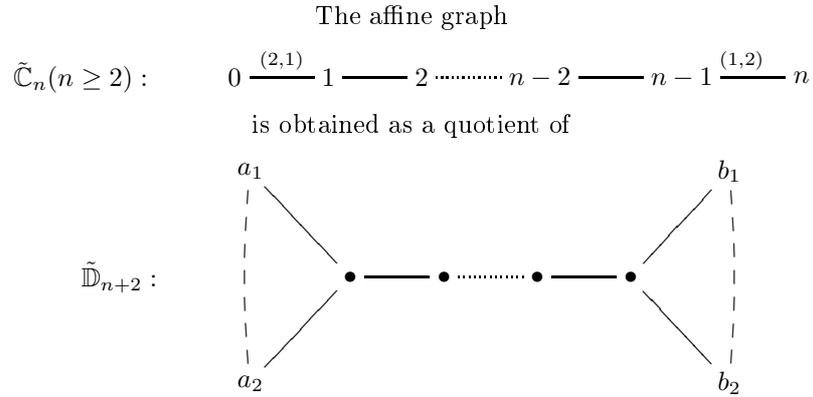
	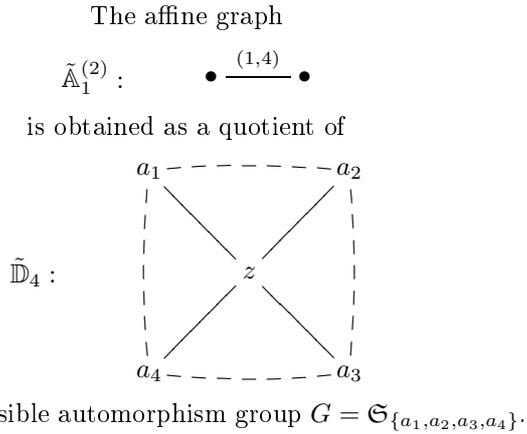
\begin{figure}[H]
		The affine graph
		$$\xymatrix{
			\tilde {\mathbb A}^{(2)}_1:  & \bullet \ar@{-}[r]^{(1,4)} & \bullet
		}$$		
		is obtained as a quotient of
		$$\xymatrix{
						& a_1 \ar@{--}@/^2pt/[rr] && a_2 \ar@{--}@/^2pt/[dd]\\
			\tilde {\mathbb D}_4:   && z \ar@{-}[rd] \ar@{-}[ld]\ar@{-}[ru]\ar@{-}[lu]\\
						& a_4 \ar@{--}@/^2pt/[uu]&& a_3 \ar@{--}@/^2pt/[ll]\\
		}$$
		by the admissible automorphism group $G=\mathfrak S_{\ens{a_1,a_2,a_3,a_4}}$. 
	\caption{Realizing $\tilde {\mathbb A}^{(2)}_1$}\label{figure:D4atoA21a}
	\end{figure}
	
	\begin{figure}[H]	
		The affine graph
		$$\xymatrix{
			\tilde{\mathbb B}_n (n\geq 2): & \bullet \ar@{-}[r]^{(1,2)} & \bullet \ar@{-}[r] & \bullet \ar@{.}[r]& \bullet \ar@{-}[r] & \bullet \ar@{-}[r]^{(2,1)} & \bullet\\
		}$$
		is obtained as a quotient of
		$$\xymatrix{
			& & a_1 \ar@{--}@/^2pt/[dd] \ar@{.}[rr] && a_{n-1} \ar@{--}@/^2pt/[dd] \ar@{-}[rd]\\ 
				\tilde {\mathbb A}_n:  & a_0=b_0 \ar@{-}[ru] \ar@{-}[rd] && & & a_n=b_n\\
			& & b_1 \ar@{.}[rr] && b_{n-1} \ar@{-}[ru]
		}$$
		by the admissible automorphism group $G=\<\prod_{i=1}^{n-1} (a_i,b_i)\>$.
	\caption{Realizing $\tilde{\mathbb B}_n$}\label{figure:AnatoBna}
	\end{figure}
	
	\begin{figure}[H]
		The affine graph 
		$$\xymatrix{
			\tilde{\mathbb C}_n (n\geq 2): & 0 \ar@{-}[r]^{(2,1)} & 1 \ar@{-}[r] & 2 \ar@{.}[r]& n-2 \ar@{-}[r] & n-1 \ar@{-}[r]^{(1,2)} & n\\
		}$$
		is obtained as a quotient of 
		$$\xymatrix{
			&a_1 \ar@{--}@/_2pt/[dd] \ar@{-}[rd]&&&&& b_1 \ar@{--}@/^2pt/[dd] \\
			\tilde{\mathbb D}_{n+2} : && \bullet \ar@{-}[r] & \bullet \ar@{.}[r] & \bullet \ar@{-}[r] & \bullet \ar@{-}[ru]\ar@{-}[rd]\\
			&a_2 \ar@{-}[ru]&&&&& b_2 \\
		}$$
		by the admissible automorphism group of automorphisms $G=\mathfrak S_{a_1,a_2} \times \mathfrak S_{b_1,b_2}$.
	\caption{Realizing $\tilde{\mathbb C}_n$}\label{figure:DnatoCna}
	\end{figure}
	
	\begin{figure}[H]
		The affine graph 
		$$\xymatrix{
			\tilde{\mathbb {BC}}_n (n\geq 2): & \bullet \ar@{-}[r]^{(1,2)} & \bullet \ar@{-}[r] & \bullet \ar@{.}[r]& \bullet \ar@{-}[r] & \bullet \ar@{-}[r]^{(1,2)} & \bullet\\
		}$$
		is obtained as a quotient of
		$$\xymatrix{
				&		&&&&&z_1 \ar@{--}@/^10pt/[dddd]\ar@{--}@/^10pt/[dddddd]\\
				&		&& b_2 \ar@{--}@/^2pt/[dddd] \ar@{.}[rr] && b_{n} \ar@{--}@/^2pt/[dddd] \ar@{-}[ru]\ar@{-}[rd]\\ 
				& 		& &&&& z_2\ar@{--}@/^10pt/[dd]\ar@{--}@/^10pt/[dddd]\\
		\tilde {\mathbb D}_{2n+2}:&  	 a\ar@{-}[rruu] \ar@{-}[rrdd]\\
				& 		& &&&& z_3\\
				&		&& c_2 \ar@{.}[rr] && c_{n} \ar@{-}[ru]\ar@{-}[rd]\\ 
				&		&&&&&z_4\\
		}$$
		by the admissible automorphism group $$G=\<\sigma(z_1,z_3)(z_2,z_4),\sigma(z_1,z_4)(z_2,z_3)\>$$
		where $\sigma=\prod_{i=2}^n (b_i,c_i)$.
	\caption{Realizing $\tilde{\mathbb {BC}}_n$}\label{figure:DnatoBCna}
	\end{figure}
	
	\begin{figure}[H]	
		The affine graph 
		$$
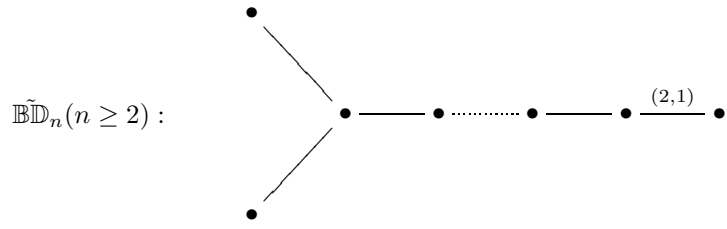
\xymatrix{
			&\bullet \ar@{-}[rd]\\
			\tilde{\mathbb {BD}}_n (n\geq 2): &&  \bullet \ar@{-}[r] & \bullet \ar@{.}[r]& \bullet \ar@{-}[r] & \bullet \ar@{-}[r]^{(2,1)} & \bullet\\
			&\bullet \ar@{-}[ru]
		}$$
		is obtained as a quotient of
		$$
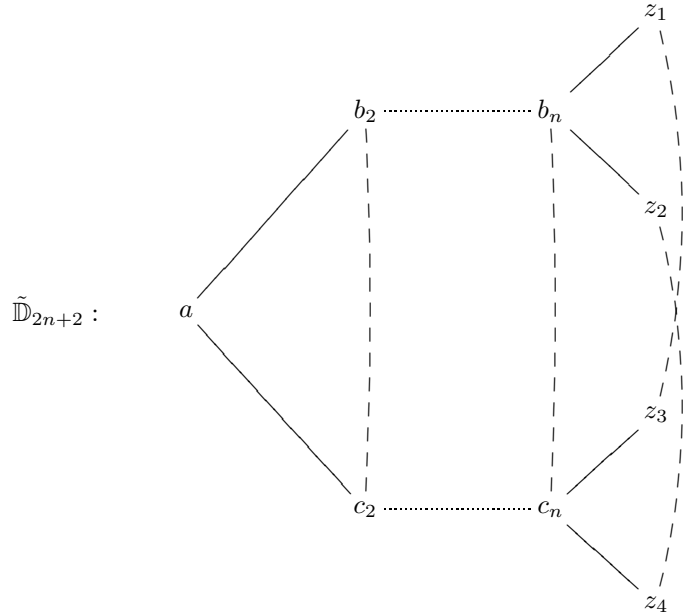
\xymatrix{
			&		&&&&&z_1 \ar@{--}@/^10pt/[dddd]\\
			&		&& b_2 \ar@{--}@/^2pt/[dddd]\ar@{.}[rr] && b_{n}\ar@{--}@/^2pt/[dddd] \ar@{-}[ru]\ar@{-}[rd]\\ 
			& 		& &&&& z_2 \ar@{--}@/^10pt/[dddd]\\
	\tilde {\mathbb D}_{2n+2}:&  	 a\ar@{-}[rruu] \ar@{-}[rrdd]\\
			& 		& &&&& z_3\\
			&		&& c_2 \ar@{.}[rr] && c_{n} \ar@{-}[ru]\ar@{-}[rd]\\ 
			&		&&&&&z_4\\
		}$$
		by the group $$G=\<\sigma(z_1,z_3)(z_2,z_4)\>$$ of admissible automorphisms where $\sigma=\prod_{i=2}^n (b_i,c_i)$.
	\caption{Realizing $\tilde{\mathbb {BD}}_n$}\label{figure:DnatoBDna}
	\end{figure}
	
	\begin{figure}[H]
		The affine graph
		$$\xymatrix{
			&\bullet \ar@{-}[rd]\\
			\tilde{\mathbb {CD}}_n : &&  \bullet \ar@{-}[r] & \bullet \ar@{.}[r]& \bullet \ar@{-}[r] & \bullet \ar@{-}[r]^{(1,2)} & \bullet\\
			&\bullet \ar@{-}[ru]
		}$$
		is obtained as a quotient of
		$$\xymatrix{
			&a_1 \ar@{-}[rd]&&&&& b_1 \ar@{--}@/^2pt/[dd]\\
			\tilde{\mathbb D}_{n+1} : && \bullet \ar@{-}[r] & \bullet \ar@{.}[r] & \bullet \ar@{-}[r] & \bullet \ar@{-}[ru]\ar@{-}[rd]\\
			&a_2 \ar@{-}[ru]&&&&& b_2 \\
		}$$
		by the admissible automorphism group $G=\<(b_1,b_2)\>$. 
	\caption{Realizing $\tilde{\mathbb {CD}}_n$}\label{figure:DnatoCDna}
	\end{figure}
	
	\begin{figure}[H]
		The affine graph
		$$\xymatrix{
			\tilde{\mathbb F}^{(1)}_{4}:  & \bullet \ar@{-}[r] & \bullet \ar@{-}[r] & \bullet \ar@{-}[r]^{(1,2)} &\bullet \ar@{-}[r] &\bullet 
		}$$
		is obtained as a quotient of 
		$$\xymatrix{
			& && c_2 \ar@{-}[d] \ar@{--}@/^2pt/[rrdd]\\
			& && c_1 \ar@{--}@/^2pt/[rd]\\
			\tilde {\mathbb E}_6: &
			a_2 \ar@{-}[r] & a_1 \ar@{-}[r] &z  \ar@{-}[r] \ar@{-}[u] &b_1\ar@{-}[r] &b_2\\
		}$$
		by the admissible automorphism group $G=\<(c_1,b_1)(c_2,b_2)\>$. 
	\caption{Realizing $\tilde{\mathbb F}^{(1)}_{4}$}\label{figure:E6atoF142a}
	\end{figure}
	
	\begin{figure}[H]
		The affine graph
		$$\xymatrix{
			\tilde{\mathbb F}^{(2)}_{4}:  & \bullet \ar@{-}[r] & \bullet \ar@{-}[r] & \bullet \ar@{-}[r]^{(2,1)} &\bullet \ar@{-}[r] &\bullet 
		}$$
		is obtained as a quotient of
		$$\xymatrix{
					&&& a_1 \ar@{--}@/^2pt/[dd] \ar@{-}[r] & a_2 \ar@{--}@/^2pt/[dd]  \ar@{-}[r]& a_3 \ar@{--}@/^2pt/[dd]\\
		\tilde {\mathbb E}_7: 	&b\ar@{-}[r] & z \ar@{-}[ru]\ar@{-}[rd]\\
					&&& c_1 \ar@{-}[r] & c_2 \ar@{-}[r]& c_3\\
		}$$
		by the admissible automorphism group $G=\<\prod_{i=1}^3(a_i,c_i)\>$. 
	
		\caption{Realizing $\tilde{\mathbb F}^{(2)}_{4}$}\label{figure:E7atoF24a}
	\end{figure}
	
	\begin{figure}[H]
		The affine graph
		$$\xymatrix{
			\tilde{\mathbb G}^{(1)}_{2}:  & \bullet \ar@{-}[r] & \bullet \ar@{-}[r]^{(1,3)} &\bullet
		}$$
		is obtained as a quotient of
		$$\xymatrix{
						& a_1 && a_2 \ar@{--}@/^2pt/[dd]\\
			\tilde {\mathbb D}_4:   && z \ar@{-}[rd] \ar@{-}[ld]\ar@{-}[ru]\ar@{-}[lu]\\
						& a_4 && a_3 \ar@{--}@/^2pt/[ll]\\
		}$$
		by the admissible automorphism group $G=\mathfrak S_{\ens{a_1,a_2,a_3}}$. 
		\caption{Realizing $\tilde{\mathbb G}^{(1)}_{2}$}\label{figure:D4atoG12a}
	\end{figure}
	
	\begin{figure}[H]
		The affine graph
		$$\xymatrix{
			\tilde{\mathbb G}^{(2)}_{2}:  & \bullet \ar@{-}[r] & \bullet \ar@{-}[r]^{(3,1)} &\bullet
		}$$
		is obtained as a quotient of
		$$\xymatrix{
			& && c_2\ar@{--}@/^2pt/[rrdd]\ar@{--}@/_2pt/[lldd] \ar@{-}[d]\\
			& && c_1\ar@{--}@/^2pt/[rd]\ar@{--}@/_2pt/[ld] \\
			\tilde {\mathbb E}_6: &
			a_2 \ar@{-}[r] & a_1 \ar@{-}[r] &z  \ar@{-}[r] \ar@{-}[u] &b_1\ar@{-}[r] &b_2\\
		}$$
		by the admissible automorphism group $G=\<(a_1,b_1,c_1),(a_2,b_2,c_2)\>$. 
		\caption{Realizing $\tilde{\mathbb G}^{(2)}_{2}$}\label{figure:E6atoG22}
	\end{figure}
	
	A corollary of theorem \ref{theorem:mutationquotient} in the context of mutation-finite graphs is:
	\begin{lem}\label{lem:quotientmutfini}
		Let $(\Delta, G)$ be a stable admissible pair such that $\Delta$ is mutation-finite. Then, $\Delta/G$ is mutation-finite.
	\end{lem}
	\begin{proof}
		It follows from theorem \ref{theorem:mutationquotient} that the map $B \fl B/G$ induces a surjective map
		$$\Mut^G(\Delta) \fl \Mut(\Delta/G).$$
		Since $\Mut^G(\Delta) \subset \Mut(\Delta)$, it follows that
		$$|\Mut(\Delta/G)| \leq |\Mut^G(\Delta)| \leq |\Mut(\Delta)| < \infty.$$
	\end{proof}
	
	\begin{rmq}
		Note that the converse of lemma \ref{lem:quotientmutfini} is false. Indeed, if $b,c \geq 3$ are integers, then $\Delta_{b,c}$ is mutation-finite whereas $K_{b,c}$ is not.
	\end{rmq}

	As an example of application, we can use this result in order to prove that valued affine graphs are mutation-finite. This gives an independent proof to a result of Seven \cite{seven:symmetrizable}.

	\begin{prop}\label{prop:affinemutfinie}
		A valued graph of Dynkin or affine type is mutation-finite.
	\end{prop}
	\begin{proof}
		Let $\Gamma$ be a valued graph of finite (resp. affine) type. Then $\Gamma$ can be written $Q/G$ where $Q$ is a quiver of finite (resp. affine) type and $G$ is an admissible automorphism of $Q$. In particular, since $Q$ is acyclic, it follows from proposition \ref{theorem:acyclicstable} that $(Q,G)$ is an admissible pair. According to theorem \ref{theorem:mutationfiniteBR} $Q$ is mutation-finite and thus, it follows from lemma \ref{lem:affinequotient} that $\Gamma=Q/G$ is also mutation-finite.
	\end{proof}
\end{section}

\section*{Acknowledgements}
	The author would like to thank Philippe Caldero, Bernhard Keller, Robert Marsh and Andrei Zelevinsky for their advices and corrections. He would also like to thank Laurent Demonet for interesting discussions and comments on the topic.


\end{document}